%% file: Negulescu_spectral.tex
\documentclass{siamart251216}

\input{ex_shared}

\ifpdf
\hypersetup{
  pdftitle={Spectral scheme for an energetic Fokker-Planck equation with  $\kappa$-distribution steady states.},
  pdfauthor={C. Negulescu, H. Parada}
}
\fi

\usepackage{comment}
\usepackage{wasysym}
\usepackage{eufrak}
\usepackage{stmaryrd}

\let\eps\varepsilon
\def\be{\begin{equation}}
\def\ee{\end{equation}}
\let\ds\displaystyle

\let\tn\textnormal

\def\RR{\mathbb{R}}

\def\NN{\mathbb{N}}

\def\cD{\mathcal D}

\def\cM{\mathcal M}

\def\cL{\mathcal L}

\def\bfv{\mathbf{v}}
\def\bfx{\mathbf{x}}

\newcommand\dd{{\,\rm d}}
\newcommand\dD{\mathrm{d}}

\begin{document}

\maketitle

% REQUIRED
\begin{abstract}
The concern of the present paper is the design of efficient numerical schemes for a specific Fokker-Planck equation describing the dynamics of energetic electrons (runaways) occurring in thermonuclear fusion plasmas. In the long-time limit, the velocity distribution function of these particles tends towards a thermal non-equilibrium $\kappa$-distribution function which is a steady-state of the considered Fokker-Planck collision operator. These $\kappa$-distribution functions have the particularity of being only algebraically decaying for large velocities, thus describing very well suprathermal populations. Our aim is to present two efficient spectral methods for the simulation of such energetic particle dynamics. The first method will be based on rational Chebyshev basis functions, rather than on Hermite basis sets, which are the basis set of choice for Maxwellian steady states. The second method is based on a different polynomial basis set, constructed via the Gram-Schmidt orthogonalisation process. These two new spectral schemes, specifically adapted to the physical context here considered, shall permit us to accurately describe the long-time asymptotics without significant numerical costs.
\end{abstract}

% REQUIRED

% REQUIRED

%%%%%%%%%%%%%%%%%%%%%%%%%%%%%%%%%%%%%%%%%%%%%%%%%%%%%%%%%%%%%%%
%%%%%%%%%%%%%%%%%%%%%%%%%%%%%%%%%%%%%%%%%%%%%%
\section{Introduction} \label{SEC1}
%%%%%%%%%%%%%%%%%%%%%%%%%%%%%%%%%%%%%%%%%%%%%%%
Modelling the dynamics of particles which are not in thermodynamic equilibrium is a very interesting and intricate problem. Plasmas encountered in astrophysics \cite{bian,liva,pierr} or fusion devices (tokamak) \cite{hase} are usually collisionless and do not attain thermal equilibrium, described by a Maxwellian velocity distribution function. Instead, such plasmas contain highly energetic particles, which are well-represented by so-called {\it "$\kappa$-distribution functions"}, {\it i.e.} power-law distributions with {\it heavy tails}. The departure from thermal equilibrium is a consequence of an out-of-equilibrium phenomenon, where the driving (diffusive) force and the drag (dissipation) force do no longer satisfy the {\it fluctuation-dissipation theorem} \cite{Kubo}, leading thus to an over-population of particles with very high kinetic energies.

The physical reason for the departure from Maxwellian equilibria is linked to the scattering cross-section dependence on the velocity, namely $\sigma(v) \sim v^{-4}$, meaning that while low-energy particles are often collisional, faster ones are generally not anymore, their slowing-down through Coulombian collisions becoming thus problematic at high speeds. An anomalous acceleration process, arising for example through turbulences or  Dreicer electric fields \cite{dre}, is then sufficient to accelerate them to infinite speeds (non-relativistic picture). In a relativistic framework, the frictional force does not really vanish for $v \rightarrow \infty$, but increases slowly with higher energies, permitting finally to bound the particle speed by the speed of light. For more physical details about all these phenomena we refer the interested reader to \cite{Gure,EC,Stahl}.\\

An accurate description of the energetic particle dynamics is of crucial importance for fusion reactor performances. In particular the description of {\it $\alpha$-particles}, which are the high energetic fusion products (Helium nuclei), is essential for self-heating reasons, whereas the description of the {\it runaway electron} dynamics is crucial for the understanding of the overall plasma stability. {\it Runaway electrons} are the main focus of this paper, and originate from the occurrence of an intense electric field during a plasma disruption, which accelerates the electrons to high energies until radiation losses can balance the acceleration of the electric field. Motivated by this rich physical context, we shall be concerned in this work with the mathematical description and numerical simulation of the runaway dynamics via a specific energetic Fokker-Planck equation.\\

\begin{figure}[ht]
\begin{center}
\includegraphics[scale=0.4]{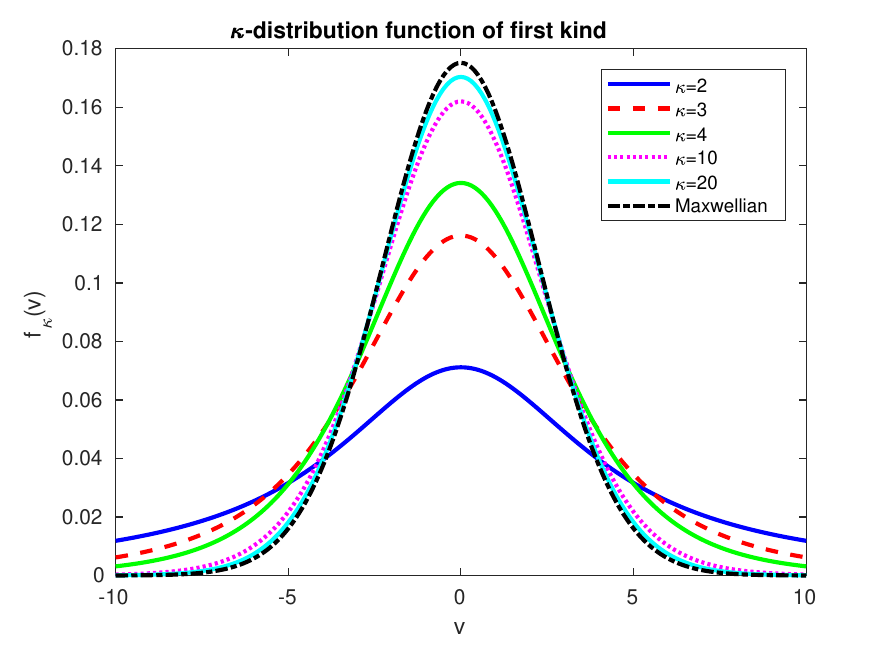}
\end{center}
\caption{$1D$ first kind $\kappa$-distribution functions for several $\kappa$-values and associated Maxwellian distribution, recovered for $\kappa \rightarrow \infty$.}\label{KM}
\end{figure}

{\it $\kappa$-distributions} were introduced for the first time in space physics by Olbert \cite{olbert} and Vasyliunas \cite{vasil} to describe suprathermal populations.
Such distributions have a core of Maxwellian type and a power-law (thus not exponential) decline in the large velocity ranges (see Fig. \ref{KM}), meaning they describe an over-population of particles at high speeds, when compared to Maxwellians. We shall consider here $\kappa$-distributions of first kind, given for $\kappa >3/2$ and $(t,\bfx,\bfv) \in (0,\infty)\times \RR^3\times \RR^3$, by the formula \cite{bian}
\be \label{kappa_first}
f_\kappa(t,\bfx,\bfv)= A_\kappa\, n(t,\bfx)\left( 1 + {|\bfv|^2 \over 2\,\kappa\, v_{th}^2}\right)^{-\kappa}\,, \qquad A_\kappa:= {1 \over (2\,\pi\, \kappa \, v_{th}^2)^{3/2}}\, {\Gamma(\kappa) \over \Gamma (\kappa-3/2)}\,, 
\ee
where the thermal speed $v_{th}:= \sqrt{k_B\, T \over m}$ and the particle density $n$ are fixed by the associated Maxwellian $\cM(t,\bfx,\bfv)$, given by
$$
\cM:={n(t,\bfx) \over (2\,\pi \, v_{th}^2)^{3/2}}\, e^{-{|\bfv|^2 \over 2\, v_{th}^2}}\,; \quad n:= \!\! \int_{\RR^3} \!\! \cM \, d \bfv = \!\! \int_{\RR^3}\!\!\!  f_\kappa \, d \bfv\,, \quad {3 \over 2}\, n\, k_B\, T := {m \over 2} \int_{\RR^3}\!\!\!  |\bfv|^2\,  \cM \, d \bfv\,.
$$
The Euler-Gamma function $\Gamma$ is defined for  $r\in (0,\infty)$ via
$$
\Gamma (r):= \int_0^\infty s^{r-1}\, e^{-s}\, ds\,, \qquad \Gamma(1)=1\,, \,\,\, \Gamma(1/2)=\sqrt{\pi}\,, \,\,\, \Gamma(r+1)=r\,\Gamma(r)\,,
$$ 
and satisfies $\Gamma(r+\alpha) \sim_{r \rightarrow \infty} \Gamma(r)\, r^\alpha$ for all $\alpha >0$, where $\sim$ stands for the asymptotic equivalence.
Remark that $f_\kappa$ contains three parameters, $v_{th}$ and $n$ which are linked to the associated Maxwellian, and $3/2 < \kappa < \infty$ which is the only free parameter and determines the distance away from the Maxwellian, in particular the lower the $\kappa$, the more pronounced are the tails. In the limit of large $\kappa \gg 1$ one recovers the standard Maxwellian distribution, namely one can easily prove that
$$
\lim_{\kappa \rightarrow \infty} f_\kappa(\bfv)=\cM(\bfv)\,, \qquad \forall \bfv \in \RR^3\,,
$$
using
$$
e^{-\xi}= \lim_{\kappa \rightarrow \infty} \left( 1 + {\xi \over \kappa}\right)^{- \kappa}\,, \qquad \lim_{\kappa \rightarrow \infty} A_\kappa = {1 \over (2\,\pi \, v_{th}^2)^{3/2}} \,.
$$

One major problem with $\kappa$-distributions is that not all the moments are bounded, contrary to Maxwellians, and this is due to the power-law decrease in the velocity-variable. Indeed, one has
$$
\int_{\RR^3} |\bfv|^\ell\, f_\kappa(\bfv) d\bfv= {2\, n \over \sqrt{\pi}} \, \left( 2\,\kappa\,v_{th}^2\right)^{\ell/2}\, { \Gamma\left({\ell+3\over 2}\right)\,\Gamma\left(\kappa-{\ell+3 \over 2}\right)\over \Gamma\left(\kappa-{3 \over 2}\right) }\,, \quad \forall\, 0 \le \ell < 2\,\kappa -3\,,
$$
higher moments being divergent, a fact which alters a characterization of the $\kappa$-distribution function at the macroscopic level (for ex. via temperature, pressure, heat flux, {\it etc.}).
To avoid this deficiency, one idea is to regularize the $\kappa$-distribution function by multiplying it with some exponentially decaying function, thus introducing $f_{\kappa,a}(\bfv):=f_\kappa(\bfv)\, e^{-a\, {|\bfv|^2 \over 2\, v_{th}^2}}$, with $a \ll 1$ a small parameter allowing to regulate the cut-off of the $\kappa$-distribution at sufficiently high speeds. The regularized $\kappa$-distribution has now a power-law behaviour at low and intermediate velocities, and an exponential cut-off at higher velocities. Besides the fact that this regularized distribution has now finite moments at all orders, this strategy is somehow also very reasonable from a physical point of view, as it can be seen as a saturation in the velocity variable, avoiding the occurrence of particles with speeds higher than the speed of light.
%Let us also remark, that in absence of the regularization or some other relativistic corrections, existing moments of the $\kappa$-distribution may also include unphysical contributions of supraluminal particles.

Let us restrict in the following to the homogeneous one-dimensional velocity case ($v \in \RR$), in which case the $\kappa$-distributions are of the form (for $\kappa >1/2$)
$$
f_\kappa(v):= c_\kappa \,n\,   \left( 1 + {v^2 \over 2 \kappa \, v_{th}^2}\right)^{-\kappa}\,, \qquad \textrm{where} \quad  c_\kappa:= {1 \over \sqrt{2\, \pi\, \kappa\, v_{th}^2 }}\, {\Gamma(\kappa) \over \Gamma (\kappa-1/2)}\,.
$$
The motivation for the restriction to the $1D$ case comes from the separation of the full $3D$ collision operator into the drag and energy scattering part (corresponding to our 1D operator) and the pitch-angle Lorentz scattering operator (skipped here), describing the scattering in the velocity direction (see \cite{NCCN} for more details).%Taking into account for the Lorentz operator requires a splitting technique, not considered here.}

%%%%%%%%%%%%%%%%%%%%%%%%%%%%%%%%%%%%%%%%%%%%%%
\subsection{Mathematical description of the energetic particle population} \label{SEC11}
%%%%%%%%%%%%%%%%%%%%%%%%%%%%%%%%%%%%%%%%%%%%%%%
Let us explain in more detail in this subsection how $\kappa$-distributions arise as equilibria of specific Fokker-Planck operators. \\
%The basic physical origin is the runaway phenomenon, which is a consequence of the fact that due to the strong energy-dependence of the Coulombian cross-section $\sigma(v)$, varying as $v^{-4}$, the dynamical friction (resistance due to long-range Coulomb interactions) of the particles cannot balance anymore the turbulent acceleration forces (or some other acceleration forces not considered here and due to the electric field), resulting in a continuous acceleration of the particles if nothing else slows them down. To translate this in mathematics, let us proceed step by step.\\

At a microscopic level, the dynamics of a particle of mass $m$ is described by Langevin's equation
$$
m\, v'(t) = - m\,\gamma(v)\, v + m\,\sqrt{2\, D(v)}\, \eta_t\,, \qquad \forall t >0\,,
$$
where $\gamma(v)=v\, \sigma(v) \sim_{|v|\gg1} v^{-3}$ is the friction coefficient and $D(v) \ge 0$ the diffusion coefficient in the velocity variable. The second term on the right-hand side stands for the random force, assumed to be a Gaussian process (white noise) described by $\eta_t$, with zero average $\langle \eta_t \rangle =0$ and infinitely short correlation times $\langle \eta_t\, \eta_{t'} \rangle =\delta_0(t-t')$. The corresponding mesoscopic picture of Langevin's equation is the kinetic Fokker-Planck equation 
$$
\partial_t f = \partial_v \left[ \gamma(v) \,v\, f + D(v)\,\partial_v f  \right]\,, \qquad \forall (t,v) \in (0,\infty) \times \RR\,.
$$
At thermal equilibrium, the friction coefficient is related to the diffusion coefficient via
${D(v) \over \gamma(v)}={k_B \, T \over m}$, which is the so-called {\it fluctuation-dissipation theorem} \cite{Kubo} and states that the friction is finally the means by which an external work is dissipated into microscopic thermal energy.
This relation leads then to the well-known form of the linear Fokker-Planck equation
\be \label{trala}
\partial_t f = \partial_v \left[ \gamma(v)\left( v\, f + v_{th}^2\,\partial_v f \right)  \right]\,, \quad \textrm{with} \quad v_{th}^2 = {k_B\, T \over m}\,, \quad \forall (t,v) \in (0,\infty) \times \RR\,.
\ee
%Remark that a transport term $v\, \partial_x f + E\, \partial_v f$ shall be introduced in this equation on the left-hand side, in order to get closer to a real thermonuclear fusion simulation with convection, for the moment we shall however only focus on the Fokker-Planck collision operator. 
Equation \eqref{trala} can be recast into the more general 1D Fokker-Planck form
\be \label{FP}
\partial_t f = \nu\, \partial_v \left[ D(v) \left( U'(v)\, f + \partial_v f \right) \right]\,, \qquad \forall (t,v) \in (0,\infty) \times \RR\,,
\ee
where $\nu >0$ is the collisional frequency and $U(v)$ a potential function describing the drift mechanism. The stationary solution of \eqref{FP} is given by
$$
f_{\infty}(v):= \alpha\, e^{ -U(v)}\,, \qquad \forall v \in \RR\,,
$$
with $\alpha >0$ determined by the initial condition $f_{\mathrm{in}}$ (conservation of mass).\\

The standard case \eqref{trala} is recovered from \eqref{FP} via $D(v)=\gamma(v)\, v_{th}^2$ and a quadratic potential $U(v):= {v^2 \over 2\, v_{th}^2}$, leading in the long-time asymptotics to the usual Maxwellian equilibrium $f_\infty(v)={n \over \sqrt{2\, \pi \, v_{th}^2}}\, e^{-{v^2 \over 2\, v_{th}^2}}$, with density $n=\int_\RR f_{\mathrm{in}}(v)\, \dD v$.\\

The fluctuation-dissipation theorem gets however violated in some situations, for example when the accelerating diffusion term cannot be anymore balanced by the friction term, leading to an out-of-equilibrium situation characterized very well by heavy-tail $\kappa$-distributions.
In mathematical terms, these $\kappa$-distributions arise by introducing in the Fokker-Planck operator, additionally to the Coulomb collisions (with a dense plasma background), an empirical turbulent acceleration mechanism described by the coefficient $\capricornus_{turb}(v)$, namely considering
\be \label{FP_astro}
\partial_t f = \nu\, \partial_v \left\{ \gamma(v)\,  \left[  v\,  f +  \left( v_{th}^2 +  \capricornus_{turb}(v) \right)\, \partial_v f \right] \right\}\,, \qquad \forall (t,v) \in (0,\infty) \times \RR\,.
\ee
The idea is that the supplementary (ad-hoc) diffusion term injects energy into the system at a rate $\capricornus_{turb}(v)$. This mechanism of energy injection, combined with the friction coefficient $\gamma(v)$, determine together the shape of the equilibrium distribution function at large velocities and permits to obtain non-equilibrium distributions.

One can now put this new Fokker-Planck equation under the form \eqref{FP} with the diffusion coefficient $D(v)= \gamma(v)\,\left(v_{th}^2  + \capricornus_{turb}(v)\right)$ and $U'(v)=  v\,\left( v_{th}^2  + \capricornus_{turb}(v)  \right)^{-1}$. For a turbulence term of the form $\capricornus_{turb}(v) =\capricornus_\star\, v^2$ with $\capricornus_\star>0$, one gets a potential of logarithmic type $U(v)=\kappa \ln \left( 1+{|\bfv|^2 \over 2\,\kappa\, v_{th}^2}\right)$, and the following $\kappa$-distribution functions as equilibria (with a constant $C_\kappa>0$)
\be \label{kappa_bis}
f_\kappa(v)= C_\kappa\, \left( 1 + {v^2 \over 2\,\kappa\,v_{th}^2}\right)^{-\kappa}\,; \qquad \kappa:= { 1  \over 2\,\capricornus_\star}\,, \quad v_{th}:= \sqrt{ k_B\, T \over m}\,.
\ee

The problem now is that this is a non-relativistic picture, the particles being able to reach infinite speeds, which is from a physical point of view not justified. To cope with this new problem, one has to introduce a saturation term, in order to restrain the speeds to speeds lower than the speed of light. From a physical point of view this arises naturally, as for very high electron speeds the friction coefficient starts to increase slowly, which is due mainly to the so-called {\it Bremsstrahlung-Effekt} \cite{Cramer,Stahl}. Starting from some threshold energy, the particles start to emit electromagnetic waves, such that their energy is dissipated, and their velocity cannot attain the speed of the light. This phenomenon can be described mathematically by starting from the previously introduced turbulent Fokker-Planck equation \eqref{FP_astro}, rewritten here as
\be \label{FP_astro_rel}
\begin{array}{lll}
%  \ds \partial_t f &=&\ds  \nu\, \partial_v \left\{ \gamma(v)\,  \left[  v\,  f +  \left( v_{th}^2 +  \capricornus_{turb}(v) \right)\, \partial_v f \right] \right\}\,, \qquad \forall (t,v) \in (0,\infty) \times \RR \\[3mm]
 \ds \partial_t f &=& \ds \nu\, \partial_v \left\{ D(v)\,  \left[  \left(v_{th}^2 +  \capricornus_{turb}(v) \right)^{-1}\, v\,  f +  \partial_v f \right] \right\}\,,
  \end{array}
\ee
and adding a new term $a\,vf$ in the friction part. This leads to
\be \label{FP_astro_rel_bis}
\partial_t f=\nu\, \partial_v \left\{ D(v)\,  \left[  \left[\left(v_{th}^2 +  \capricornus_{turb}(v) \right)^{-1}+a \right]\, v\,  f +  \partial_v f \right] \right\}\,,
\ee
with $D(v)\!=\!\gamma(v)\left(v_{th}^2  + \capricornus_{turb}(v)\right)$. 
Thus $U'(v)\!=\! v\left[\left(v_{th}^2 +  \capricornus_{turb}(v) \right)^{-1}+a \right]$, and taking again $\capricornus_{turb}(v) =\capricornus_\star\, v^2$ leads to $U(v)= \kappa \ln \left( 1+{|\bfv|^2 \over 2\,\kappa\, v_{th}^2}\right) + a\, {v^2 \over 2}$ and equilibria
$$
f_{\kappa,a}(v):=N_\kappa\, \left( 1 + {v^2 \over 2\,\kappa\,v_{th}^2}\right)^{-\kappa}\, e^{-a\, {v^2 \over 2}}\,, \qquad \forall v \in \RR\,,
$$
with a constant $N_\kappa >0$. Equation \eqref{FP_astro_rel_bis} can be rewritten via these equilibria as
\be \label{eq_D}
\partial_t f=\nu\, \partial_v \left\{ D(v)\, f_{\kappa,a}(v) \partial_v \left( {f\over f_{\kappa,a}} \right) \right\}\,, \qquad \forall v \in \RR\,.
\ee

\begin{figure}[ht]
\begin{center}
\includegraphics[scale=0.4]{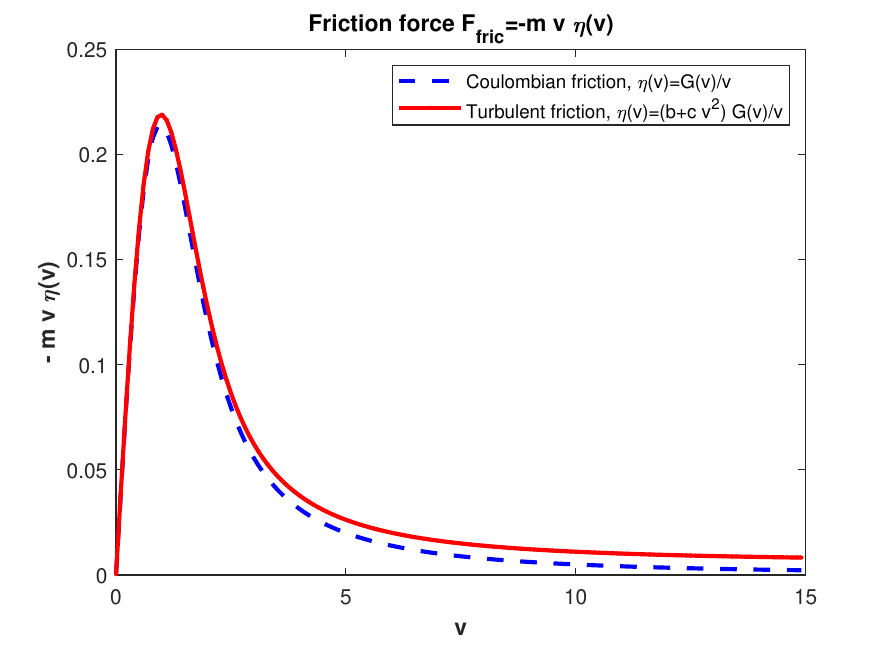}
\end{center}
\caption{Friction $F_{fric}= - m\,v \, \gamma(v)\,\left[1+a\,\left(v_{th}^2 +  \capricornus_{turb}(v) \right) \right]$ with $\gamma(v)=G(v)/v$. Coulombian friction for $a=0$ (dotted curve, vanishing at infinity); turbulent friction for $a=0.012$ (full curve).}\label{fric}
\end{figure}

The new friction force reads now $F_{fric}=-m\,\gamma(v)\, v\, \left[1+a\,\left(v_{th}^2 +  \capricornus_{turb}(v) \right) \right]$, which means that the friction does no longer vanish  for $|v| \rightarrow \infty$, but saturates at a constant value, as illustrated in Fig. \ref{fric}, where we recall that the Coulombian friction satisfies $\gamma(v) \sim {G(v) \over v} \sim_{|v|\gg 1} v^{-3}$ \cite{CHENF,Ruther}, whereas the turbulent diffusion has the form $\capricornus_{turb}(v) =\capricornus_\star\, v^2$. Here $G$ is the Chandrasekhar function \cite{EC}, defined via the error function $\phi$ by $G(v):= {\phi(v)-v\, \phi'(v) \over 2\, v^2}$ for $v \neq 0$ and $G(0)=0$.
%%%%%%%%%%%%%%%%%%%%%%%%%%%%%%%%%%%%%%%%%%%%%%
\subsection{Aim of this paper and outline} \label{SEC12}
%%%%%%%%%%%%%%%%%%%%%%%%%%%%%%%%%%%%%%%%%%%%%%%
The main objective of this paper is  the efficient numerical resolution of the following $1D$ evolution problem
	\begin{align}\label{ev_gen}
		\partial_t f = \partial_v\left[f_{\mathrm{eq}}(v)\,\partial_v\left(\frac{f(t,v)}{f_{\mathrm{eq}}(v)}\right)\right]\,,
	\qquad \forall (t,v)\in (0,\infty)\times\RR\,,\end{align}
        associated with some initial condition $f_{\mathrm{in}}$. Three cases shall be considered and compared, associated to the following three equilibria:
        \begin{itemize}
        \item[(i)] Standard Fokker-Planck case: $f_{\mathrm{eq}}=\cM(v)$;
        \item[(ii)] Energetic Fokker-Planck case: $f_{\mathrm{eq}}=f_\kappa(v)$;
        \item[(iii)] Regularized energetic Fokker-Planck case: $f_{\mathrm{eq}}=f_{\kappa,a}(v):=f_\kappa(v)\, e^{-a\, {v^2 \over 2}}$, where $a \ll 1$ is a small cut-off parameter.
        \end{itemize}
        Let us remark here that, contrary to the physical part, we set $\nu \equiv 1$, $v_{th} \equiv 1$ and $D(v) \equiv 1$ for simplicity reasons (compare with \eqref{eq_D}), the main difficulties being still present in this simplified framework. The mathematical study of the general Fokker-Planck equation has been investigated in \cite{EC}. The Maxwellian equilibrium reads
	\be \label{nr_M}
		\cM(v):=\frac{1}{\sqrt{2\,\pi}}\,e^{-v^2/2}\,, \qquad \forall v\in\RR\,,
	\ee
whereas the $\kappa$-equilibrium is given for $\kappa > 1/2$ and $ v \in \RR$ by the formula
\be \label{nr_k}
f_\kappa(v):= {c_\kappa \over  \left( 1 + {v^2 \over 2 \kappa }\right)^{\kappa}}\,, \qquad \textrm{where} \quad  c_\kappa:= {1 \over \sqrt{2\, \pi\, \kappa }}\, {\Gamma(\kappa) \over \Gamma (\kappa-1/2)}\,,
\ee
such that $\int_\RR f_\kappa (v)\, \dD v =1$. The parameter $0<a\ll 1$ shall be chosen so to restrain the particle velocities to speeds lower than the speed of the light. Remark that
$$
f_{\kappa,a} \,\, \longrightarrow_{a \rightarrow 0} \,\, f_\kappa \longrightarrow_{\kappa \rightarrow \infty} \cM\,,
$$
and we recall, for later use, the moments of the $\kappa$-distribution function
$$
\int_{\RR} |v|^\ell\, f_\kappa(v) \dd v= {1 \over \sqrt{\pi}} \, \left( 2\,\kappa\right)^{\ell/2}\, { \Gamma\left({\ell+1\over 2}\right)\,\Gamma\left(\kappa-{\ell+1 \over 2}\right)\over \Gamma\left(\kappa-{1 \over 2}\right) }\,, \qquad \forall\, 0 \le \ell < 2\,\kappa -1\,,
$$
as well as those of the regularized $\kappa$-distribution, this time for all $\ell \in \NN$
\be \label{whitt}
%\begin{array}{l}
  \ds \int_{\RR} |v|^\ell\, f_{\kappa,a}(v) \dd v= 
  %\ds {2^{\ell+1 \over 2}\, a^{-{\ell+1\over 2}} \over \Gamma(\kappa)}\, E(\kappa;b;{\ell+1 \over 2};b;a\, \kappa)\\[3mm]
 % \hspace{1cm}=
  \ds 2^{\ell+1 \over 2}\, \Gamma\left({\ell+1\over 2}\right)\,a^{{\kappa \over 2}-{\ell+1\over 4}-{1 \over 2}} \, \kappa^{{\kappa \over 2}+{\ell+1\over 4}-{1 \over 2}}\, e^{a\, \kappa \over 2}\, W_{{1 \over 2}-{\kappa \over 2}-{\ell+1\over 4},{\kappa \over 2}-{\ell+1\over 4}}(a\, \kappa) \,,
%\end{array}
\ee
where $W(\cdot)$ is the Whittaker function \cite{Zw}.\\
%$E(\cdot)$ is the Mac Robert's E-function with $b \in \RR$ an arbitrary value, and $W(\cdot)$ is the Whittaker function \cite{Zw}.\\

The first immediately visible difficulty in a numerical resolution of \eqref{ev_gen} is related to the unboundedness of the velocity-domain, which leads to challenging numerical complications. Several possibilities are available to treat such unbounded problems, namely spectral methods, domain truncation, mapping techniques (see \cite{boyd} for more details). We shall concentrate in the present paper on spectral methods, based on basis functions intrinsic to unbounded domains, such as for example Hermite-basis functions, Sinc-functions, rational Chebyshev functions, {\it etc}. Spectral methods are so-called {\it global methods} and, when well-designed, are superior in performance to local methods, such as standard finite-difference or finite-element methods, due to their exponential accuracy (greater precision with smaller number of points or modes). This exponential accuracy is however only achieved when the basis functions are well-chosen, an inadequate choice can have a drastic impact on the convergence rate.

The mathematical and numerical treatment of the three problems (i)-(iii) in view of the construction of a spectral scheme is completely different, and this is mainly due to the fact that the spectrum of the standard case (i) is discrete, the eigenvectors being well-known (Hermite basis), whereas case (ii) admits a continuous spectrum, requiring a continuous spectral representation with generalized eigenvectors not belonging to the functional space of interest. A basis not constituted of eigenvectors is then of rescue in this case (ii). The equilibrium $\kappa$-distributions are however not exponentially decaying in the velocity variable, but rather algebraically, such that the usual exponentially decaying Hermite polynomials (useful in case (i)) are not suitable here. When a function $f(t,v)$ decays as an inverse power of the velocity-variable, then rational Chebyshev polynomials are clearly more adequate, such that our spectral method for case (ii) is based on the introduction of a rational Chebyshev basis set.

Finally, case (iii) admits again (due to the cut-off in the high velocity ranges) a discrete spectrum, however with non-standard eigenvectors. But the regularization allows now to employ the function $f_{\kappa,a}$ as a valid weight function, to start the construction (via Gram-Schmidt orthogonalisation) of an adequate polynomial basis set, permitting by this manner to introduce an efficient spectral method for case (iii).\\

The outline of this paper is the following. In Section \ref{SEC2} we state the mathematical problem and present some of its properties (spectrum, asymptotic decay). Section \ref{SEC3} is concerned with the construction of a spectral scheme for the resolution of the Fokker-Planck equation in case (ii), scheme based on rational Chebyshev polynomials (so-called {\it RC-scheme}). Section \ref{SEC4} focuses on the regularized Fokker-Planck equation of case (iii) and proposes a spectral method based on a polynomial basis set, constructed via the weight function $f_{\kappa,a}$ (so-called {\it GS-scheme}). Finally, Section \ref{SEC5} presents the numerical results obtained with both spectral schemes and compares them with a standard finite-difference method ({\it FD-scheme}).

%%%%%%%%%%%%%%%%%%%%%%%%%%%%%%%%%%%%%%%%%%%%%%
\section{Fokker-Planck collision operator} \label{SEC2}
%%%%%%%%%%%%%%%%%%%%%%%%%%%%%%%%%%%%%%%%%%%%%%%
In this section we recall briefly some mathematical results from \cite{EC}, permitting to understand the construction of our spectral numerical schemes and their validation. We consider the following $1D$ evolution problem
\be\label{evolution_pb_generic}
(FP)\,\,\,  \left\{
  \begin{array}{l}
    \ds \partial_t f = \partial_v\left[f_{\mathrm{eq}}(v)\,\partial_v\left(\frac{f}{f_{\mathrm{eq}}}\right)\right]\,,\qquad \forall (t,v)\in (0,\infty)\times\RR\,,\\[3mm]
    \ds f(0,\cdot)=f_{\mathrm{in}}(\cdot)\,,
  \end{array}
  \right.
  \ee
where for the equilibrium distribution function we choose either $f_{\mathrm{eq}}(v)=\cM(v)$ defined in \eqref{nr_M}, $f_{\mathrm{eq}}(v)=f_\kappa(v)$ defined in \eqref{nr_k} or the regularized version $f_{\mathrm{eq}}(v)=f_{\kappa,a}(v)$.  In the long-time limit $t \rightarrow \infty$ the solution to this equation will tend towards the stationary solution given by $f_\infty(v):={\langle f_{\mathrm{in}} \rangle \over \langle f_{\mathrm{eq}} \rangle}\, f_{\mathrm{eq}}(v)$, where the average in the velocity space will be denoted by
$$
\langle f \rangle:=\int_\RR f(v)\dd v\,.
$$
The second order differential operator occurring in \eqref{evolution_pb_generic}, {\it i.e.}
\be \label{op_L}
\cL_{\mathrm{eq}}(f):= -\partial_v\left[f_{\mathrm{eq}}\,\partial_v\,\left(\frac{f}{f_{\mathrm{eq}}}\right)\right]\,, \qquad  \cL_{\mathrm{eq}}: \cD(\cL_{\mathrm{eq}}) \subset L^{2}_{f^{-1}_{\mathrm{eq}}} \rightarrow L^{2}_{f^{-1}_{\mathrm{eq}}}\,,
\ee
is a {\it linear, unbounded, self-adjoint and positive operator} on the Hilbert-space
		\begin{align} \label{HL}
			L^{2}_{f^{-1}_{\mathrm{eq}}}:=\left\{f:\RR \rightarrow \RR\,\,\,,\quad \int_\RR \vert f\vert^2\,f^{-1}_{\mathrm{eq}}\,\dd v<\infty\right\}, \qquad (f,g)_{f^{-1}_{\mathrm{eq}}}:=\int_\RR fg\,f^{-1}_{\mathrm{eq}}\,\dd v\,,
		\end{align}
and with domain 
		\begin{align}
			\cD(\cL_{\mathrm{eq}})=\{f\in L^2_{f^{-1}_{\mathrm{eq}}},\quad \cL_{\mathrm{eq}}(f)\in L^2_{f^{-1}_{\mathrm{eq}}}\}\,.
		\end{align}
              
Standard functional analysis (Lion's theorem, Hille-Yosida theorem) then permits showing that there exists a unique solution  $f \in C^1([0,\infty);L^2_{f^{-1}_{\mathrm{eq}}}) \cap C([0,\infty);\cD(\cL_{\mathrm{eq}}))$ to the Fokker-Planck equation \eqref{evolution_pb_generic}, for each initial condition $f_{\mathrm{in}} \in \cD(\cL_{\mathrm{eq}})$.

Before going on, let us also remark that one has the inclusions
$$
L^2_{\cM^{-1}} \subset L^2_{f_{\kappa,a}^{-1}} \subset L^2_{f_{\kappa}^{-1}} \subset L^2(\RR)\,.
$$
    
	%%%%%%%%%%%%%%%%%%%%%%%%%%%%%%%%%%%%%%%%%%%%%%
	\subsection{Standard Fokker-Planck equation with $f_{\mathrm{eq}}=\cM$} \label{SEC21}
	%%%%%%%%%%%%%%%%%%%%%%%%%%%%%%%%%%%%%%%%%%%%%%%
	This case is very well documented \cite{Grad,anto}, we recall the results here only for completeness reasons. The equation we consider for $(t,v)\in (0,\infty)\times\RR$ is simply
        \be \label{St_FP}
(FP)_\cM\,\,\,          \left\{
  \begin{array}{l}
\ds \partial_t f = -\cL_\cM(f)=\partial_v\left[\cM\,\partial_v\left(\frac{f}{\cM}\right)\right]= \partial_v\left[ v\, f + \partial_v f\right] \,,\\[3mm]
\ds f(0,\cdot)=f_{\mathrm{in}}(\cdot)\,.
  \end{array}
  \right.
        \ee
	%%%%%%%%%%%%%%%%%%%%%%%%%%%%%%%%%%%%%%%%%%
	%%%%%%%%%%%%%%%%%%%%%%%%%%%%%%%%%%%%%%%%%%%
	\begin{proposition} \cite{EC,anto} {\bf (Hermite basis functions)} \label{PropH}
		The operator $\cL_\cM$ defined in \eqref{op_L} is self-adjoint, positive and with compact resolvent, such that its spectrum is discrete, real, positive and consists of a sequence of eigenvalues $(\lambda_k)_{k\in\NN} \subset \RR$ satisfying $\lambda_k\to\infty$ as $k\to\infty$. In particular one has for all $k \in \NN$
$$
\cL_{\cM}\, \psi_k(v) = \lambda_k\, \psi_k(v)\,, \qquad \forall v \in \RR\,,   \quad \textrm{with eigenvalue}\,\,\, \lambda_k:=k\,, 
$$
and the associated eigenvectors are the Hermite functions, defined recursively via $\psi_0\equiv \cM\,,\,\, \psi_1\equiv v\,\cM$ and for all $k\ge 1$ by the formulae
		\begin{align}
			\sqrt{k+1}\, \psi_{k+1}(v)=v\,\psi_k(v)-\sqrt{k}\,\psi_{k-1}(v)\,, \quad \forall v \in \RR\,.
		\end{align}
Remark that $\{\psi_k\}_{k\in\NN}$ form a complete orthonormal basis set of the space $L^2_{\cM^{-1}}$.
	\end{proposition}
	%%%%%%%%%%%%%%%%%%%%%%%%%%%%%%%%%%%%%%%%%%
	%%%%%%%%%%%%%%%%%%%%%%%%%%%%%%%%%%%%%%%%%%%
	This proposition is useful from an analytic point of view, permitting to study the decay rate towards the equilibrium, but also from a numerical perspective, permitting the construction of a spectral method. 

\begin{proposition} \cite{EC,anto} {\bf (Time decay)} \label{PropD}
Let $f$ be the unique solution of the evolution problem \eqref{St_FP}, with initial condition $f_{\mathrm{in}}  \in L^2_{\cM^{-1}}$. Then the following exponential time-decay of the solution towards its stationary state holds
	\begin{align}\label{theoretical_decay_case1}
		\Vert f(t)-f_{\infty}\Vert_{\cM^{-1}}\le\Vert f_{\mathrm{in}}-\langle f_{\mathrm{in}} \rangle\, \cM\Vert_{\cM^{-1}} \,e^{-t}\,, \quad \forall t \ge 0\,,
	\end{align}
        where the stationary solution is given by $f_\infty(v)=\langle f_{\mathrm{in}} \rangle\, \cM(v)$.
\end{proposition}
	
Now, one can make use of the orthonormal basis of eigenvectors $\{\psi_k\}_{k\in\NN} \subset L^2_{\cM^{-1}}$  to expand the solution of the evolution problem \eqref{St_FP}, as
$$
		f(t,v)=\sum_{k=0}^\infty  \alpha_k(t)\,\psi_k(v)\,,\qquad 
		\alpha_k(t)=(f(t),\psi_k)_{\cM^{-1}}=\int_\RR f(t,v)\,\psi_k(v)\,\cM^{-1}(v)\dd v \,.
        $$
Inserting this expression in \eqref{St_FP} yields an equation to be solved for the expansion coefficients $\alpha_k(t)$, leading for all $t \ge 0$ and $k \in \NN$ to 
$$
\alpha_k'(t) + \lambda_k\, \alpha_k(t)=0 \quad \Rightarrow \quad \alpha_k(t)=e^{-kt}\,\alpha_{\tn{in},k}\,, \quad \alpha_{\tn{in},k}:=\int_\RR f_{\mathrm{in}}(v)\,\psi_k(v)\,\cM^{-1}(v)\dd v\,.
$$
This spectral decomposition is the starting point of spectral methods for standard Fokker-\-Planck equations. We refer the interested reader to \cite{gi,Wang21} for further details.
%%%%%%%%%%%%%%%%%%%%%%%%%%%%%%%%%%%%%%%%%%%
\subsection{Energetic Fokker-Planck equation for  $f_{\mathrm{eq}}=f_\kappa$} \label{SEC22}
%%%%%%%%%%%%%%%%%%%%%%%%%%%%%%%%%%%%%%%%%%%
We are now coming to the investigation of the energetic Fokker-Planck equation, which writes 
\be \label{EEQQ}
(FP)_\kappa\,\,\,\left\{
  \begin{array}{l}
    \ds  \partial_t f = -\cL_\kappa(f)=\partial_v\left[f_\kappa(v)\,\partial_v\left(\frac{f}{f_\kappa(v)}\right)\right]= \partial_v\left[  {v \over 1+{v^2 \over 2\, \kappa}}\, f + \partial_v f \right]\,, \\[3mm]
    \ds f(0,\cdot)=f_{\mathrm{in}}(\cdot)\,,
  \end{array}
  \right.
\ee
where the equilibrium is this time $f_\kappa(v)= c_\kappa\, (1+{v^2 \over 2\,\kappa})^{-\kappa}$, with $c_\kappa >0$ the normalization constant. We shall work in the Hilbert-space $L^2_{f_\kappa^{-1}}$ for some fixed $\kappa \in (0,\infty)$.\\

The linear, unbounded operator $\cL_\kappa$, defined in \eqref{op_L},
%\be \label{Kappa_FP}
%\cL_\kappa(f):=-\partial_v\left[f_\kappa(v)\,\partial_v\left(\frac{f}{f_\kappa}\right)\right]\,, \qquad \cL_\kappa: \cD(\cL_\kappa) \subset L^{2}_{f^{-1}_\kappa} \rightarrow L^{2}_{f^{-1}_\kappa}\,,
%        \ee
has no more a discrete spectrum, in contrast to  the standard Fokker-Planck case (see Proposition \ref{PropH}). Nevertheless, one can now follow the same strategy as in the standard Fokker-Planck case, by searching for {\it generalized eigenfunctions} and a continuous spectral representation, as has been done in \cite{EC}. In the present paper however, we shall follow a different strategy, by considering the expansion of the $\kappa$-distribution function in a basis set, not constituted of (generalized) eigenfunctions of the corresponding Fokker-Planck operator $\cL_\kappa$, but constituted of a different basis set of the functional Hilbert-space $L^2_{f_\kappa^{-1}}$, namely of rational Chebyshev polynomials introduced in Section \ref{SEC3}. For the moment, let us finish this Section, by recalling a time-decay result of the solution to \eqref{EEQQ} towards its stationary solution. Due to the fact that there is no spectral gap, the time-decay rate will be no more exponential as in the standard Fokker-Planck case, but algebraic.
	
	\begin{theorem}\label{Theorem_Dolbeault}\cite{Dol_Bouin_Ziv,EC} {\bf (Time decay)}
		Let $f$ be a solution of problem \eqref{EEQQ}, with initial condition $f_{\mathrm{in}}\, f_\kappa^{-1}\in L^\infty(\RR_v)$. Then for all $0<p<2\,\kappa-1$, with $\kappa > 1$, the following estimate holds
		\begin{align}\label{convergence_dolbeault}
			\Vert f(t)- f_\infty \Vert_{f^{-1}_\kappa}^2 \le \left[\Vert f_{\mathrm{in}}-\langle f_{\mathrm{in}} \rangle\, f_\kappa  \Vert_{f^{-1}_\kappa}^{-4/p} + K_{p,\kappa}\, {4 t \over p} \right]^{-p/2} \,, \qquad \forall t \ge 0\,,
		\end{align}
with $f_\infty=\langle f \rangle\, f_\kappa=\langle f_{\mathrm{in}} \rangle\, f_\kappa$ and $K_{p,\kappa}>0$ a constant depending on $p$ and $\kappa$.
	\end{theorem}
        %%%%%%%%%%%%%%%%%%%%%%%%%%%%%%%%%%%%%%%%%%%	
	\subsection{Regularized energetic Fokker-Planck operator for $f_{\mathrm{eq}}=f_{\kappa,a}$}\label{SEC23}
        %%%%%%%%%%%%%%%%%%%%%%%%%%%%%%%%%%%%%%%%%%%
        Finally, let us turn to the investigation of the regularized Fokker-Planck equation
\be \label{EEQQ_bb}
(FP)_{\kappa,a}\,\,\,\left\{
  \begin{array}{l}
    \ds  \partial_t f = - \partial_v\left[f_{\kappa,a}(v)\,\partial_v\left(\frac{f}{f_{\kappa,a}}\right)\right]= \partial_v\left[  \left( {1 \over 1+{v^2 \over 2\, \kappa}} + a \right)\, v\, f + \partial_v f \right]\,,\\[3mm]
    \ds f(0,\cdot)=f_{\mathrm{in}}(\cdot)\,,
  \end{array}
  \right.
\ee
with equilibrium $f_{\kappa,a}(v)= c_\kappa\, (1+{v^2 \over 2\,\kappa})^{-\kappa}\, e^{-a{v^2 \over 2}}$, $c_\kappa >0$ being the normalization constant of $f_\kappa$. We shall work in $L^2_{f_{\kappa,a}^{-1}}$, with some fixed $\kappa \in (0,\infty)$ and $0<a \ll 1$.\\

To study the spectrum of this Fokker-Planck equation, the Liouville transformation $f= g\, \sqrt{f_{\kappa,a}}$ is very useful, permitting to switch from the eigenvalue problem $\cL_{\kappa,a}(f) = \lambda\, f$ to the Schr\"odinger problem (sharing the same spectrum)
$$
- \partial_{vv} g + Q_a(v)\, g = \lambda \, g\,, \qquad \forall v \in \RR\,,
$$
where the occurring potential $Q_a$ is a confining potential given by
$$
Q_a(v):= {v^2 \over 4} \left[ {1 \over 1+{v^2 \over 2\, \kappa}} + a \right]^2 - {1 \over 2}\, \left[ {1-{v^2 \over 2\, \kappa} \over (1+{v^2 \over 2\, \kappa})^2} + a \right] \sim_{|v| \gg 1} {1 \over 4}\, a^2\, v^2 - {a \over 2}\,.
$$
Spectral theorems of quantum mechanics \cite{Pearson1988} permit now to show that the spectrum of this Schr\"odinger equation is discrete for $a>0$ (it is continuous for $a=0$), and this essentially due to confining properties of the potential, {\it i.e.} $Q_a(v) \rightarrow_{|v| \rightarrow \infty} \infty$. We recall that for $a=0$ one has $Q_0(v) \sim_{|v| \gg 1} {\kappa\, (\kappa+1) \over v^2}$ \cite{EC}, hence not confining. The main goal of the regularization is thus to have a discrete spectrum. Furthermore, the Fokker-Planck operator $\cL_{\kappa,a}$ being self-adjoint and positive, one has $\sigma(\cL_{\kappa,a})=\{\lambda_k\}_{k \in \NN} \subset [0,\infty)$, with $\lambda_0=0$ the first eigenvalue, corresponding to the eigenfunction $f_{\kappa,a}$ of $\cL_{\kappa,a}$.\\

The existence of a spectral gap of the operator $\cL_{\kappa,a}$ is strictly related to the existence of a Poincar\'e inequality, leading to an exponential time-decay behaviour.
\begin{comment}
To investigate now the decay in time of the solution to \eqref{EEQQ_bb} towards the stationary state $f_\infty= {\langle f_{\mathrm{in}} \rangle \over \langle f_{\kappa,a} \rangle}\,f_{\kappa,a} $, it is more suitable to filter out the equilibrium by setting $f=h\,f_{\kappa,a}$ and to consider the so-called {\it Ornstein-Uhlenbeck} evolution equation for the unknown $h$ in the weighted Hilbert-space $L^2_{f_{\kappa,a}}$
\be \label{EEQQ_hh}
(OU)_{\kappa,a}\,\,\,\left\{
  \begin{array}{l}
    \ds  \partial_t h = {1 \over f_{\kappa,a}}\,  \partial_v\left[f_{\kappa,a}(v)\,\partial_v h \right]\,, \quad \forall (t,v) \in (0,\infty) \times \RR\,,\\[3mm]
    \ds h(0,\cdot)=h_{\mathrm{in}}(\cdot)=f_{\mathrm{in}}(\cdot)\,f^{-1}_{\kappa,a}\,.
  \end{array}
  \right.
\ee
This equation has the same spectrum as the original Fokker-Planck equation $(FP)_{\kappa,a}$, exhibiting thus a spectral gap which leads to a Poincar\'e inequality of the type
$$
||h-{\overline h}||^2_{L^2_{f_{\kappa,a}}} \le C_P ||\partial_v h||^2_{L^2_{f_{\kappa,a}}}\,, \qquad \forall h \in H^1_{f_{\kappa,a}}\,,
$$
with $C_P>0$ the Poincar\'e constant and the average defined as ${\overline h}:= \int_\RR h\, f_{\kappa,a}\, \dD v$. Hence, one has the following exponential time-decay property.
\end{comment}
\begin{theorem}\label{Thm_kappa}{\bf (Time decay)}
  Let $f$ be a solution of problem \eqref{EEQQ_bb}, with initial condition $f_{\mathrm{in}}\in L^2_{f_{\kappa,a}^{-1}}$. Then the following estimate holds, for some constant $C>0$
  \begin{align}\label{theoretical_decay_case1_kappa}
		\Vert f(t)-f_{\infty}\Vert_{f_{\kappa,a}^{-1}}\le C\, \Vert f_{\mathrm{in}}-{\langle f_{\mathrm{in}} \rangle \over \langle f_{\kappa,a} \rangle}\, f_{\kappa,a} \Vert_{f_{\kappa,a}^{-1}} \,e^{-\lambda_{1,a}\,t}\,, \quad \forall t \ge 0\,,
	\end{align}
        where the stationary solution is given by $f_\infty= {\langle f_{\mathrm{in}} \rangle \over \langle f_{\kappa,a} \rangle}\,f_{\kappa,a} $ and $\lambda_{1,a}>0$ is the first positive eigenvalue of $\cL_{\kappa,a}$, satisfying $\lambda_{1,a} \longrightarrow_{a \rightarrow 0} 0$.
  \end{theorem}
        %%%%%%%%%%%%%%%%%%%%%%%%%%%%%%%%%%%%%%%%%%%
%%%%%%%%%%%%%%%%%%%%%%%%%%%%
\section{Spectral scheme based on rational Chebyshev functions} \label{SEC3}
%%%%%%%%%%%%%%%%%%%%%%%%%%%%%
We are coming now to the design of an adequate spectral method for our energetic particle Fokker-Planck equation \eqref{EEQQ}, we recall here for clarity reasons
\be \label{EEQQ_again}
(FP)_\kappa\,\,\,\left\{
  \begin{array}{l}
    \ds  \partial_t f = \partial_v\left[f_\kappa(v)\,\partial_v\left(\frac{f}{f_\kappa(v)}\right)\right]= \partial_v\left[  {v \over 1+{v^2 \over 2\, \kappa}}\, f + \partial_v f \right]\,,\\[3mm]
    \ds f(0,\cdot)=f_{\mathrm{in}}(\cdot)\,.
  \end{array}
  \right.
\ee
In the previous section we evoked the Hermite spectral method for the resolution of the standard Fokker-Planck equation \eqref{St_FP}. Hermite functions however are exponentially decaying for large $v$, and hence their use is not natural for the approximation of functions which are slowly decaying for large $v$, such as the $\kappa$-equilibrium of our problem \eqref{EEQQ_again}. New basis sets are thus required and have to be introduced. We shall not focus on basis functions composed of generalized eigenfunctions (continuous spectral representation) as has been done in \cite{EC}. Using instead transformations that map a finite domain into an infinite domain, allows to generate a broad class of new basis functions for the unbounded domain, that are images of, for example, Chebyshev polynomials. We shall focus on this strategy in this section. 
%%%%%%%%%%%%%%%%%%%%%%%%%%%%%%%%%%%%%%%%%%%	
	\subsection{Rational Chebyshev polynomials}\label{SEC31}
%%%%%%%%%%%%%%%%%%%%%%%%%%%%%%%%%%%%%%%%%%%	
        At the ground of our polynomial basis construction are the Chebyshev polynomials, which are no more than a change-of-variable in the Fourier cosine functions. Indeed, the map $s \in [0,\pi] \mapsto \xi=\cos(s) \in [-1,1]$ transforms the cosine functions into the {\it Chebyshev polynomials} $T_n(\xi)$, {\it i.e.}
\be \label{Cheb}
T_n(\xi):= \cos(n\,s) \quad \textrm{where} \quad  s=\arccos(\xi)\,, \qquad \forall \xi \in [-1,1]\,.
\ee
The thus defined set $\{T_n\}_{n \in \NN}$ forms an orthogonal basis of $L^2((-1,1); \rho(\xi)d\xi)$ with the weight $\rho(\xi):= { 1 \over \sqrt{1-\xi^2}}$, in particular one has the orthogonality relations
        $$
        \int_{-1}^1 {T_n(\xi)\, T_m(\xi) \over \sqrt{1-\xi^2}}\, d\xi=
        \left\{
        \begin{array}{lcl}
          \ds 0\,, & \textrm{if} & \ds m \neq n\\[1mm]
          \ds \pi\,, & \textrm{if} & \ds m=n=0\,, \qquad \qquad \forall n,m \in \NN\,.\\[1mm]
          \ds {\pi \over 2}\,, & \textrm{if} & \ds m=n \neq 0\,,
        \end{array}
        \right.
        $$
%        A Chebyshev series is finally just a Fourier cosine expansion with a change-of-variable, in order to map a periodic interval $[0,\pi]$ into a bounded, non-periodic interval $[-1,1]$.
        
%        Let us at the end also remark that the Chebyshev polynomials are the eigenfunctions of the singular Sturm-Liouville problem
%        $$
%-\sqrt{1-\xi^2}\, \partial_\xi \left[\sqrt{1-\xi^2} \, \partial_\xi T_n(\xi) \right] = n^2\, T_n(\xi)\,, \quad \forall \xi \in (-1,1)\,.
%        $$

\noindent At this point, using now the algebraic map
\be \label{algM_1}
\xi \in (-1,1) \mapsto v:= {L\, \xi \over \sqrt{1-\xi^2}} \in \RR\,, \quad \textrm{or equiv.} \quad  \xi:= {v \over \sqrt{L^2+v^2}}\,,
\ee
 or
\be \label{algM_2}
s \in (0,\pi) \mapsto v:= L\, \cot(s)\in \RR\,, \quad \textrm{or equiv.} \quad s=\textrm{arccot}(v/L)\,,
\ee
permits to transform the finite interval into an infinite interval, and one defines the {\it rational Chebyshev polynomials} as
\be \label{TB}
{\mathcal C}^L_n(v):=T_n(\xi)=\cos(n\, s)\quad \textrm{where} \quad s \in (0,\pi)\,,\,\,\, \xi \in (-1,1)\,, \,\,\,v \in \RR\,. 
\ee
One can equally start from the sinus functions in order to define
\be \label{SB}
{\mathcal S}^L_n(v):=S_n(\xi)=\sin((n+1)\, s)\quad \textrm{where} \quad s \in (0,\pi)\,,\,\,\, \xi \in (-1,1)\,, \,\,\,v \in \RR\,. 
\ee
Let us remark that $L>0$ is a mapping parameter, which has to be adjusted as closely as possible to the scale of the function to be expanded.\\

The two sets $\{{\mathcal C}^L_n\}_{n \in \NN}$ along with $\{{\mathcal S}^L_n\}_{n \in \NN}$ form two independent orthogonal basis sets of the Hilbert-space $L^2_{\sigma_L}$, with weight $\sigma_L(v):= {L\over L^2+v^2}$. In particular one has
$$\int_{-\infty}^\infty {\mathcal C}^L_n(v)\, {\mathcal C}^L_m(v) {L \over L^2+v^2}\, \dD v=\left\{
\begin{array}{lcl}
\ds 0\,, & \textrm{if} & \ds m \neq n\\[1mm]
\ds \pi\,, & \textrm{if} & \ds m=n=0\,, \qquad \qquad \forall n,m \in \NN\,,\\[1mm]
\ds {\pi \over 2}\,, & \textrm{if} & \ds m=n \neq 0\,,
\end{array}\right.$$
and equally for $\{{\mathcal S}^L_n\}_{n \in \NN}$. For an efficient computation of these rational Chebyshev basis functions, one has at hand the following recurrence formulae
        $$
{2\, v/L \over \sqrt{1+(v/L)^2}}\, {\mathcal C}^L_n(v)= {\mathcal C}^L_{n+1}+{\mathcal C}^L_{n-1}\,, \quad ({\mathcal C}^L_n)'(v)={n \over 2\, L}\, {1 \over \sqrt{1+(v/L)^2}} \left[ {\mathcal C}^L_{n-1}-{\mathcal C}^L_{n+1}\right]\,,
$$
as well as
$$
{4\, (v/L)^2 \over 1+(v/L)^2}\, {\mathcal C}^L_n(v)= {\mathcal C}^L_{n+2}+2\, {\mathcal C}^L_n+{\mathcal C}^L_{n-2}\,, \quad {4 \over 1+(v/L)^2}\, {\mathcal C}^L_n(v)= -{\mathcal C}^L_{n+2}+2\, {\mathcal C}^L_n-{\mathcal C}^L_{n-2}\,.
$$

Using all these formulae permits to compute the first polynomials, which read
$$
{\mathcal C}^L_0 \equiv 1\,, \quad {\mathcal C}^L_1(v)={{v/L} \over \sqrt{1+(v/L)^2}}\,, \quad {\mathcal C}^L_2={(v/L)^2 -1 \over 1+(v/L)^2}\,, \quad {\mathcal C}^L_3= { (v/L)\, [(v/L)^2-3] \over [1+(v/L)^2]^{3/2}}\,,
$$
and
$$
{\mathcal S}^L_0(v)={1 \over \sqrt{1+(v/L)^2}}\,, \quad {\mathcal S}^L_1(v)= {2\,(v/L) \over 1+(v/L)^2}\,, \quad {\mathcal S}^L_2(v)= {3(v/L)^2 - 1 \over [1+(v/L)^2]^{3/2}}\,.
$$
As an example, consider the simple rational (algebraically decaying) functions $[1+(v/L)^2]^{-n}$ as well as  $[1+(v/L)^2]^{-(n+1/2)}$ which are rather harmless, however approximating them with Hermite polynomials fails drastically. In this case rational Chebyshev polynomials are more powerful, as can be observed through the following relations, valid for all $n \in \NN$

        \be \label{expi}
        { 1 \over [1+(v/L)^2]^n}= \alpha_0\, {\mathcal C}^L_0+  \cdots + \alpha_{2\, n} \, {\mathcal C}^L_{2\, n}\,, \quad { 1 \over [1+(v/L)^2]^{n+1/2}}= \beta_0\, {\mathcal S}^L_0+ \cdots + \beta_{2\,n} \, {\mathcal S}^L_{2\,n}\,.
        \ee
        The next theorem describes the way of convergence of the rational Chebyshev series.
        \begin{definition} An expansion of the form
          $$
f(v)= \sum_{n=0}^\infty \alpha_n\, \varphi_n(v)\,, \qquad \forall v \in \RR\,,
$$
is said to be algebraically or exponentially convergent, if there exist some constants $p>0$, $\mu >0$, $\iota>1$ and $N_\star \in \NN$, such that one of the two items is satisfied
\begin{itemize}
\item $|\alpha_n|\le p\, n^{-\iota}\,, \quad \forall n \ge N_\star \quad \textrm{(algebraic convergence)}\,,$
  \item $|\alpha_n|\le p\, e^{-\mu\, n}\,, \quad \forall n \ge N_\star \quad \textrm{(exponential convergence)}\,.$
  \end{itemize}
          \end{definition}

    \vspace{0.2cm}
    
        \begin{theorem} \cite{boyd} \label{bboo}{\bf (Convergence of the rational Chebyshev expansion)}\\
          Let $g \in L^2_{\sigma_L}$ be a function free of singularities on $\RR$. Then one has:\\
          (i) If $g$ has an asymptotic power series for $|v| \gg 1$ containing only even, non-negative powers of $1/v$, the expansion of $g$ in the basis $\{{\mathcal C}^L_n(v)\}_{n \in \NN}$ is exponentially convergent.\\
          (ii) If $g$ has an asymptotic power series for $|v| \gg 1$ that contains only odd, positive powers of $1/v$, the expansion of $g$ in the basis $\{{\mathcal S}^L_n(v)\}_{n \in \NN}$ is exponentially convergent.\\
          (iii) If $g$ has an asymptotic power series for $|v| \gg 1$ that contains only integral, non-negative powers of $1/v$, the expansion of $g$ in the full basis $\{{\mathcal C}^L_n(v),{\mathcal S}^L_n(v)\}_{n \in \NN}$ is exponentially convergent.
          %(iv) If $g$ has an asymptotic power series for $|v| \gg 1$ that contains non-integral powers of $1/v$, the expansion series of $g$ in one of the basis $\{{\mathcal C}^L_n(v)\}_{n \in \NN}$ or/and $\{{\mathcal S}^L_n(v)\}_{n \in \NN}$ is only algebraically convergent.
          \end{theorem}

%%%%%%%%%%%%%%%%%%%%%%%%%%%%%%%%%%%%%%%%%%%	
	\subsection{RC-spectral method}\label{SEC32}
        %%%%%%%%%%%%%%%%%%%%%%%%%%%%%%%%%%%%%%%%%%%
	We shall set now $L= \sqrt{2\, \kappa}$ and take $\kappa \in \NN$. The question is to find a basis set for our weighted Hilbert-space $L^2_{f_\kappa^{-1}}$, based on the just introduced rational Chebyshev basis sets of the space $L^2_{\sigma_\kappa}$, where we recall that 
        $$
f_\kappa(v)={ c_\kappa \over (1+{v^2 \over 2\, \kappa})^{\kappa}}\,, \qquad c_\kappa:= {1 \over \sqrt{2\, \pi\, \kappa }}\, {\Gamma(\kappa) \over \Gamma (\kappa-1/2)}\,, \qquad \sigma_\kappa(v):={1 \over\sqrt{2\, \kappa}}\, {1\over 1+{v^2 \over 2\, \kappa}}\,.
$$
Scaling the rational Chebyshev polynomials $\{{\mathcal C}^L_n, \, {\mathcal S}^L_n \}_{n \in \NN}$ as 
$$
\Theta_n(v):= {\mathcal C}^L_n(v)\, \Upsilon_\kappa(v)\,, \qquad \Xi_n(v):= {\mathcal S}^L_n(v)\, \Upsilon_\kappa(v)\,, 
$$
with
$$
\Upsilon_\kappa(v):= \left[ \sigma_\kappa(v)\, f_\kappa(v)\right]^{1/2}={(c_\kappa/\sqrt{2\, \kappa})^{1/2} \over (1+{v^2 \over 2\, \kappa})^{\kappa+1 \over 2}}\,, \qquad \forall v \in \RR\,,
$$
permits to introduce two independent orthogonal basis sets $\{\Theta_n, \, \Xi_n \}_{n \in \NN}$ for $L^2_{f_\kappa^{-1}}$. The recurrence formulae remain the same as above, the only modification is in the derivative, which becomes
$$
\Theta'_n(v)={n \over 2}\, {1 \over\sqrt{2\, \kappa}}\,{1 \over \sqrt{1+{v^2 \over 2\, \kappa}}} \left[ \Theta_{n-1}(v)-\Theta_{n+1}(v)\right] -  {\kappa +1 \over 2\, \kappa}\, {v \over 1 + {v^2 \over 2\, \kappa}}\, \Theta_n(v)\,, \quad \forall v \in \RR\,.
$$
Taking an odd $\kappa \in \NN$ prescribes an expansion in $\{\Theta_n\}_{n \in \NN}$ (see Thm. \ref{bboo} (i)) and symmetry properties permit furthermore to show that  even functions $f(t,|v|)$ possess an expansion using solely even rational basis functions. Altogether we shall thus search for a solution to \eqref{EEQQ_again} under the form
$$
f(t,v)= \sum_{n=0}^\infty \alpha_{2n}(t)\, \Theta_{2n}(v)\,, \qquad \forall (t,v) \in [0,\infty)\times \RR\,.
$$

Introducing this spectral expansion in the FP-equation \eqref{EEQQ_again}, permits to find an ODE system for the spectral coefficients $\{\alpha_{2n}(t)\}_{n \in \NN}$. Straight computations yield 
\begin{equation}
\label{eq: LTB}
\begin{array}{lll}
-{\mathcal L_\kappa}(\Theta_0)&\!\!\!=&\!\!\!c_{0}\, \Theta_{0}+d_{0}\, \Theta_{2}+e_{0}\, \Theta_{4}\,,\\[3mm]
%-{\mathcal L_\kappa}({\mathcal C}^L_1)=&c_{1}\, {\mathcal C}^L_{1}+d_{1}\, {\mathcal C}^L_{3}+e_{1}\, {\mathcal C}^L_{5}\,,\\
-{\mathcal L_\kappa}(\Theta_2)&\!\!\!=&\!\!\!b_{2}\, \Theta_{0}+c_{2}\, \Theta_{2}+d_{2}\, \Theta_{4}+e_{2}\, \Theta_{6}\,,\\[3mm]
%-{\mathcal L}({\mathcal C}^L_3)=&b_{3}\, {\mathcal C}^L_{1}+c_{3}\, {\mathcal C}^L_{3}+d_{3}\, {\mathcal C}^L_{5}+e_{3}\, {\mathcal C}^L_{7}\,,\\[3mm]
-{\mathcal L_\kappa}(\Theta_{2n})\!\!&\!\!\!=&\!\!\!\!a_{2n}\, \Theta_{2n-4}+b_{2n}\, \Theta_{2n-2}+c_{2n}\, \Theta_{2n}+d_{2n}\, \Theta_{2n+2}+e_{2n}\, \Theta_{2n+4}\,, \,\,\, \forall n\geq2\,,
\end{array}
\end{equation}\

with 
\[c_{0}=-{(\kappa-1)^2\over 16\, \kappa}\,, \quad d_{0}=-{\kappa-1 \over 4\, \kappa}\,, \quad e_0={(\kappa-1)\, (\kappa+3) \over 16\, \kappa},\]
%\[c_{1}={-6\kappa-3\over8}\,, \quad d_{1}={9-22\kappa\over 16}\,, \quad e_1={-3+2\kappa\over16},\]
\[b_{2}=-{\kappa-1\over 8\, \kappa}\,, \quad c_{2}=-{25 + \kappa\,(\kappa-6)\over 32\, \kappa}\,, \quad d_{2}=-{\kappa-9 \over 8\, \kappa}\,, \quad e_2={(\kappa-3)\,(\kappa+5) \over 32\, \kappa},\]
%\[b_{3}={2\kappa+21\over16}\,, \quad c_{2}={4\kappa-27\over8}\,, \quad d_{2}={12-5\kappa\over 4}\,, \quad e_2={10\kappa-15\over16},\]
and
\[a_{n}={ (\kappa-1+n)\,(\kappa+3-n) \over 32\, \kappa}\,, \quad b_{n}=-{\kappa -  (n-1)^2 \over 8\, \kappa}\,, \quad c_n={-3\,n^2 - (\kappa-1)^2 \over 16\, \kappa}\]
$$
d_{n}=-{\kappa - (n+1)^2 \over 8\, \kappa}\,, \quad e_{n}={ (\kappa-1-n)\,(\kappa+3+n) \over 32\, \kappa}\,, \quad \forall n\geq4\,.
$$
Taking then the weighted scalar-product (in $L^2_{f_\kappa^{-1}}$) of the Fokker-Planck equation $\partial_t f = -{\mathcal L_\kappa} f$ with $\Theta_{2j}$, yields

\[\begin{cases}
 \alpha_{0}^{\prime}(t)=&\!\!\!\!\!c_{0}\alpha_{0}(t)+b_{2}\alpha_{2}(t)+a_{4}\alpha_{4}(t), \\
 &\vdots\\
 \alpha_{2j}^{\prime}(t)=&\!\!\!\!\!e_{2j-4}\alpha_{2j-4}(t)+d_{2j-2}\alpha_{2j-2}(t)+c_{2j}\alpha_{2j}(t)+b_{2j+2}\alpha_{2j+2}(t)+a_{2j+4}\,\alpha_{2j+4}\,.
\end{cases}\]

\vspace{0.2cm}

This leads to the following penta-diagonal ODE system to be solved for the computation of the spectral coefficients ${\mathcal X}(t):= (\alpha_0(t),\alpha_2(t), \alpha_4(t), \cdots)^t$, {\it i.e.}
\be \label{L_syst}
{\mathcal X}'(t)= - M\, {\mathcal X}(t)\,, 
\ee
with
{\small $$
M= -
\left(
\begin{array}{cccccccc}
  c_0&b_2&a_4&0&\cdots&0&0&0\\[2mm]
  d_0&c_2&b_4&a_6&\cdots&0&0&0\\[2mm]
  e_0&d_2&c_4&b_6&\cdots&0&0&0\\[2mm]
  \vdots&&&\ddots&\ddots&&&\vdots\\[2mm]
  0&0&0&&\cdots&d_{2N-4}&c_{2N-2}&b_{2N}\\[2mm]
  0&0&0&&\cdots&e_{2N-4}&d_{2N-2}&c_{2N}\\[2mm]
\end{array}
\right)\,.
$$}
\begin{remark} \label{remCheb}
The just introduced spectral method is very  accurate for $\kappa$-values arising in astrophysics or fusion plasmas, such as $\kappa \in (1,10)$, however it turns out to be problematic when $\kappa$ becomes larger and larger, preventing the investigation of the limit $\kappa \rightarrow \infty$. Indeed, one can immediately observe that some entries of the matrix $M$ become infinite for $\kappa \rightarrow \infty$ and the rational basis functions $\{{\mathcal C}^L_n\}_{n \in \NN}$ get also more and more similar in this limit, rendering the RC-spectral method not adequate to treat cases with $\kappa \gg 1$. Furthermore, the expansion of the asymptotic equilibrium $f_\kappa(v)$ requires the use of $\kappa$ Chebyshev basis terms, owing to
$$
f_\kappa(v)=\sum_{n=0}^N \alpha_{2n}\, {\mathcal C}_{2n}(v)\, \Upsilon_\kappa(v) \,\,\, \Leftrightarrow \,\,\, f_\kappa(v)\, \Upsilon^{-1}_\kappa(v)= { \sqrt{c_\kappa}\, (2\, \kappa)^{1/4} \over (1+{v^2 \over 2\, \kappa})^{\kappa-1 \over 2}}=\sum_{n=0}^N \alpha_{2n}\, {\mathcal C}_{2n}(v)\,,
$$
such that with \eqref{expi} one notices that this is possible for $2 N=\kappa-1$. Thus, with $\kappa \rightarrow \infty$ the expansion requires more and more terms, leading to an inefficient Chebyshev-spectral method. There are however some advantages for small $\kappa$ values, discussed in the numerical section, such that we find that this scheme has its right to exist.
\end{remark}
%%%%%%%%%%%%%%%%%%%%%%%%%%%%
\section{Spectral scheme based on non-standard orthogonal polynomials} \label{SEC4}
%%%%%%%%%%%%%%%%%%%%%%%%%%%%%
The aim of this section will be thus to construct a different spectral scheme, which shall degenerate into the Hermite spectral method of Section \ref{SEC21}, when $\kappa$ tends towards infinity. The idea of this new scheme is based on the Gram-Schmidt orthogonalisation procedure, which requires the existence of the moments of the weight function ($\kappa$-distribution), such that we shall be forced to consider in this section rather the following regularized Fokker-Planck equation for all $(t,v) \in (0,\infty) \times \RR$
\be \label{Kappareg_FP}
(FP)_{\kappa,a}\,\,\,\left\{
  \begin{array}{l}
    \ds  \!\! \partial_t f = - \cL_{\kappa,a}(f)= \partial_v \left[\left(\dfrac{1}{1+\frac{v^{2}}{2\kappa}}+a\right)\,vf+\partial_{v}f\right]=  \partial_v\left[ \, f_{\kappa,a}\, \partial_v \left( f/f_{\kappa,a} \right) \right]\\[3mm]
    \ds \!\! f(0,\cdot)=f_{\mathrm{in}}(\cdot)\,.
  \end{array}
  \right.
\ee
Even if the operator $\cL_{\kappa,a}$ has a discrete spectrum as the standard Fokker-Planck one, the eigenfunctions are non-standard and have no analytical expression, such that one has to search for them numerically. We will follow however a distinct strategy, constructing a different orthogonal basis $\{\zeta_k(v)\}_{k\in\mathbb{N}}$  of $L^{2}_{f^{-1}_{\kappa,a}}$, based on Gram-Schmidt's orthogonalization process, allowing the expansion of the solution to  \eqref{Kappareg_FP} as
\be \label{exp_GS}
f(t,v)= \sum_{k=0}^\infty \alpha_k(t)\, \zeta_k(v)\,, \qquad \forall (t,v) \in [0,\infty) \times \RR\,,
\ee
with $\{\alpha_k(t)\}_{k \in \NN}$ to be determined.\\
%The main idea is to construct a spectral method which shall turn into the Hermite spectral method of section \ref{SEC21}, when $\kappa$ tends towards infinity.\\

To simplify the presentation, let us filter out the equilibrium from the solution to \eqref{Kappareg_FP} by setting $f=h\, f_{\kappa,a}$, where the unknown $h \in L^2_{f_{\kappa,a}}$ is this time solution of the so-called {\it Ornstein-Uhlenbeck} equation
\be \label{Kappareg_FP_h}
(OU)_{\kappa,a}\,\,\,\left\{
  \begin{array}{l}
    \ds  \partial_t h ={1 \over  f_{\kappa,a} }\, \partial_v \left[\, f_{\kappa,a} \, \partial_v h  \right] \,, \qquad \forall (t,v) \in (0,\infty) \times \RR\,,\\[3mm]
    \ds h(0,\cdot)=h_{\mathrm{in}}(\cdot)=f_{\mathrm{in}}\,f^{-1}_{\kappa,a}\,.
  \end{array}
  \right.
\ee
We shall construct now an orthogonal polynomial basis set $\{p_k(v)\}_{k\in\mathbb{N}}$  of $L^{2}_{f_{\kappa,a}}$, considering $\omega(v):=f_{\kappa,a}(v)$ as weight function. The regularization is crucial here, as this construction requires the existence of all moments of the weight function.
The  Gram-Schmidt orthogonalization process yields the following recursive formula for the construction of a monic polynomial set
\be \label{rec1}
\begin{gathered}
p_{-1}\equiv 0\,, \qquad p_{0}\equiv1\,, \qquad  p_1(v)=v\,,\\[3mm]
p_{k+1}(v)=vp_{k}(v)-\beta_{k}p_{k-1}(v), \qquad \forall k\geq 1\,,
\end{gathered}
\ee
where 
$$
\beta_k:= { (v\, p_k, p_{k-1})_\omega \over || p_{k-1}||^2_\omega}={ || p_{k}||^2_\omega\over || p_{k-1}||^2_\omega}\,, \qquad \forall k\geq 1\,.
$$
Classical weights, such as Maxwellian ones, give rise to classical families of orthogonal polynomials (Hermite polynomials for ex.). For such families explicit formulae are available for the recursion coefficients. Outside the classical framework however, numerical strategies are required for the computation of these coefficients, and the difficult part in the construction of the polynomial basis set is the stable computation of these coefficients $\{\beta_k\}_{k\in\mathbb{N}}$. One has to take into account that we deal with polynomials on unbounded domains, and small (round-off or truncation)  errors are immediately visible. In particular, polynomial spectral methods are more sensitive to errors than for ex. Fourier spectral methods. We shall present a stable computational strategy for the recursion coefficients below, for the moment however, let us remark that the first coefficients write
\be \label{Han}
\beta_{1}=\dfrac{m_{2}}{m_{0}}\,, \qquad \beta_{2}=\dfrac{m_{0}m_{4}-m_{2}^{2}}{m_{0}m_{2}}\,; \qquad m_\ell:= \int_\RR v^\ell\, \omega(v)\, \dD v\,, \qquad \forall \ell \in \NN\,,
\ee
where $\{m_\ell\}_{\ell \in \NN}$ are the moments of order $\ell$ of the weight function $\omega(v)$, computed via Whittaker formula \eqref{whitt}. Supposing now that we have constructed in this manner an  orthogonal basis set $\{p_k\}_{k\in\mathbb{N}}$  of $L^{2}_{f_{\kappa,a}}$, the solution to \eqref{Kappareg_FP_h} writes
\be \label{rec}
h(t,v)= \sum_{k=0}^\infty \alpha_k(t)\, p_k(v)\,, \qquad \forall (t,v) \in [0,\infty) \times \RR\,.
\ee
Inserting this decomposition into Ornstein-Uhlenbeck's evolution equation \eqref{Kappareg_FP_h} and taking the $L^2_\omega$ scalar product with $p_\ell(v)$, permits to get a system for the computation of the spectral coefficients $\{\alpha_k(t)\}_{k\in\mathbb{N}}$. Indeed, one gets
$$
\sum_{k=0}^\infty \alpha_k'(t) \int_{\RR} p_k(v)\, p_\ell(v)\, \omega(v)\, \dD v = - \sum_{k=0}^\infty \alpha_k(t) \int_{\RR} p'_k(v)\, p'_\ell(v)\, \omega(v)\, \dD v\,, \qquad \forall \ell \in \NN\,.
$$
Remarking that $(p_k,p_\ell)_\omega= \gamma_\ell\, \delta_{k,\ell}$, where we denote the weighted norm by $\gamma_\ell:=||p_\ell ||^2_\omega=\int_\RR |p_\ell(v)|^2\, \omega(v)\, \dD v$\,,
one obtains the linear ODE system
{\small
\be \label{ODE}
\left( 
\begin{array}{c}
\alpha'_0(t)\, \gamma_0\\[2mm]
\alpha'_1(t)\, \gamma_1\\[2mm]
\vdots \\[2mm]
\alpha'_\ell(t)\, \gamma_\ell\\[2mm]
\vdots
\end{array}
\right)= -
\left( 
\begin{array}{cccc}
\theta_{00}&\cdots&\theta_{0k}&\cdots\\[2mm]
\vdots&\ddots&\vdots& \\[2mm]
\theta_{\ell 0}&\cdots&\theta_{\ell k}&\cdots\\[2mm]
\vdots&\vdots&\vdots&\ddots
\end{array}
\right)\, 
\left( 
\begin{array}{c}
\alpha_0(t)\\[2mm]
\alpha_1(t)\\[2mm]
\vdots \\[2mm]
\alpha_k(t)\\[2mm]
\vdots
\end{array}
\right)\,,
\ee}
where the entries of the matrix are
\be \label{INTeg}
\theta_{k\ell}=\theta_{\ell k}:=\int_{\RR}  p'_k(v)\, p'_\ell(v)\, \omega(v)\, \dD v\,.
\ee
It remains now to find an accurate and stable manner to compute the integrals $\{\theta_{\ell k}\}_{\ell,k\in\mathbb{N}}$, the recurrence coefficients $\{\beta_k\}_{k\in\mathbb{N}}$ as well as the weighted norms $\{\gamma_k\}_{k\in\mathbb{N}}$.\\

Several methods have been proposed in literature for an accurate computation of the recursion coefficients, see for example \cite{Gautschi}. The classical approach in terms of moments (see \eqref{Han}), computed via the Hankel determinants, is numerically problematic, being very ill-conditioned. In this work we  shall use the so-called \textit{Modified Chebyshev algorithm}, which is a more stable and a better conditioned strategy. The idea is to work with "modified moments", corresponding to a different measure (weight function) $d \mu:=\mu(v)\, \dD v$, to be chosen close to the original measure $d \omega:=\omega(v)\, \dD v$, and with the additional requirement of possessing explicit recurrence coefficients. This shall permit to bypass the stability issues and generate in a stable manner a non-classical orthogonal polynomial basis set.

To be more precise, taking for instance as weight $\mu(v)=\cM(v)\, e^{-a\,\frac{ v^{2}}{2}}$,  which is close to our original weight $\omega_\kappa(v)= f_\kappa\, e^{-a\,\frac{ v^{2}}{2}}$ in the sense that
%\be \label{mu_ome}
$\lim_{ \kappa \rightarrow \infty} \omega_\kappa = \mu$,
%\ee
we construct the corresponding orthogonal monic polynomial basis set $\{q_k\}_{k\in\mathbb{N}}$ of $L^2_\mu(\RR)$, via the three-term recurrence relation
\be \label{rec2}
\begin{gathered}
q_{-1} \equiv 0, \qquad q_{0}\equiv 1\,, \qquad q_1(v)=v\,, \\[2mm]
q_{k+1}(v)=vq_{k}(v)-b_{k}q_{k-1}(v), \qquad k\geq 1\,.
\end{gathered}
\ee
For such Maxwellian weights $\mu(v)$, the recursion coefficients $\{b_{k}\}_{k \in \NN}$ are explicitly known and $\{q_{k}\}_{k \in \NN}$ are Hermite-like polynomials, in particular one has
$$
b_k={k \over 1+a}\,, \qquad ||q_k||_\mu^2={k! \over (1+a)^{2k+1 \over 2}}, \qquad k\geq 1\,.
$$
Based on this new basis set, let us introduce the \textit{mixed moments}
\[\sigma_{k,\ell}=\int_{\mathbb{R}}p_{k}(v)q_{\ell}(v)\, \omega(v)\,\dD v, \qquad \forall k,\ell\geq-1\,,\]
which will help to construct the original basis set $\{p_{k}\}_{k \in \NN}$.
By orthogonality, we have $\sigma_{k,\ell}=0$ for $k>\ell$. Moreover, as $p_k$ and $q_k$ are monic polynomials of degree $k$, we get
\[p_k(v)=q_{k}(v)+r(v), \qquad \text{ with deg}(r)<k,\]
such that using the orthogonality of $\{p_{k}\}_{k \in \NN}$, one has
\be \label{gam}
\gamma_k=\|p_{k}\|^{2}_{\omega}=\int_{\mathbb{R}}p_{k}(v)\,p_k(v)\, \omega(v)\, \dD v=\int_{\mathbb{R}}p_{k}(v)\,q_{k}(v)\, \omega(v)\,\dD v=\sigma_{k,k}.
\ee
Thus
\be \label{bet}
\beta_{k}={||p_k||^2_\omega \over ||p_{k-1}||^2_\omega }=\dfrac{\sigma_{k,k}}{\sigma_{k-1,k-1}}\,, \quad \forall k \ge 1\,.
\ee
Using the recurrence formulae \eqref{rec1} for $p_k$ and \eqref{rec2} for $q_k$, we derive
\be \label{sig}
\sigma_{k,\ell}=\sigma_{k-1,\ell+1}-\beta_{k-1}\,\sigma_{k-2,\ell}+b_{\ell}\,\sigma_{k-1,\ell-1}\,,  \quad \forall k,\ell \ge 1\,.
\ee
Starting hence from the values of $\sigma_{-1,\ell} =0$ and $\{\sigma_{0,\ell}\}_{\ell \in \NN}= \int_{\mathbb{R}} q_{\ell}(v)\, \omega(v)\,\dD v$, which are nothing more than linear combinations of the moments $\{m_{\ell}\}_{\ell \in \NN}$, computed via the Whittaker formulae \eqref{whitt}, one can compute step by step the values of $\{\sigma_{k,\ell}\}_{k,\ell \in \NN}$, $\{\beta_k\}_{k \in \NN}$ and $\{\gamma_{k}\}_{k \in \NN}$ via \eqref{gam}-\eqref{sig}.\\

For the resolution of our Fokker-Planck problem \eqref{Kappareg_FP_h}, thus of \eqref{ODE}, one needs now to compute the entries of the matrix, in particular $\theta_{k\ell}:=\int_{\RR}  p'_k(v)\, p'_\ell(v)\, \omega(v)\, \dD v$. This shall be done by introducing the following quantities
$$
\chi_{k,\ell}:= \int_{\RR}  p_k(v)\, p'_\ell(v)\, \omega(v)\, \dD v\,, \quad \xi_{k,\ell}:= \int_{\RR}  p'_k(v)\, p_\ell(v)\, \omega(v)\, \dD v\,.
$$
Using now the recursion formulae \eqref{rec1}-\eqref{rec2} permits to obtain for all $k,\ell \in \NN$
$$
\begin{array}{lll}
\theta_{k,\ell}&=&\ds \theta_{k-1,\ell+1}+ \beta_\ell\, \theta_{k-1,\ell-1}+\chi_{k-1,\ell}- \beta_{k-1}\,\theta_{k-2,\ell}-\xi_{k-1,\ell}\,,\\[3mm]
\chi_{k,\ell}&=&\ds \chi_{k-1,\ell+1}+ \beta_\ell\, \chi_{k-1,\ell-1}-\gamma_\ell\,\delta_{k-1,\ell}- \beta_{k-1}\,\chi_{k-2,\ell}\,, \\[3mm]
\xi_{k,\ell}&=& \ds \xi_{k-1,\ell+1}+ \beta_\ell\, \xi_{k-1,\ell-1}+\gamma_\ell\, \delta_{k-1,\ell}- \beta_{k-1}\,\xi_{k-2,\ell}\,.
\end{array}
$$
Thus starting from the values $\theta_{0,\ell}=0$, $\chi_{-1,\ell}=0$, $\xi_{0,\ell}=0$ for all $\ell \in \NN$, as well as from the values one has to compute as an initial step
$$
\theta_{1,\ell}:= \int_{\RR}   p'_\ell(v)\, \omega(v)\, \dD v\,, \quad \chi_{0,\ell}:= \int_{\RR}   p'_\ell(v)\, \omega(v)\, \dD v\,, \quad \xi_{1,\ell}:= \int_{\RR}  p_\ell(v)\, \omega(v)\, \dD v\,,
$$
one gets the desired values of $\{\theta_{k,\ell}\}_{k,\ell \in \NN}$. The computation of $\theta_{1,\ell}$ shall be done via Whittaker's formulae \eqref{whitt}, for accuracy reasons.
Remark also that $\xi_{1,\ell}=\gamma_0\, \delta_{0,\ell}$ due to the orthogonality of the polynomials $\{p_k\}_{k \in \NN}$, and that $\chi_{k,\ell}=\xi_{\ell,k}$ for all $k,\ell \in \NN$. Furthermore, due to
$$
p'_\ell(v)=\sum_{k=0}^{\ell-1} c_k^{(\ell)}\, p_k(v)\,, \quad \textrm{with coefficients} \,\,\, c_k^{(\ell)}= {\xi_{\ell,k} \over \gamma_k}\,,
$$
one has immediately that $\chi_{k,\ell} =0$ for all $k > \ell-1$ as well as $\xi_{k,\ell}=0$ for all $\ell > k-1$. Finally, this decomposition permits also to observe that
$$
\chi_{0,\ell}=c_0^{(\ell)} \,\gamma_0 =\theta_{1,\ell}\,, \qquad c_0^{(\ell)}:=(p'_\ell,p_0)_\omega\,.
$$
\begin{remark} \label{remGS}
Let us mention that the Gram-Schmidt spectral method is very performant for large $\kappa$-values, $\kappa \gg 1$,  as in this case the weight $\omega_\kappa$ is close to the Maxwellian weight $\mu$, whereas for small values of $\kappa$ the use of Whittaker's formulae \eqref{whitt} becomes problematic due to the fact that we deal with polynomials on unbounded domains, rendering the utility of the {\it GS}-scheme very awkward. There is a delicate interplay between the regularization parameter $a$ (thus somehow the truncation of the domain) and the $\kappa$-parameter (meaning the hot tails of the distribution function) one has to care of, in order to avoid instability problems. Let us also observe that the asymptotic equilibrium $f_{\kappa,a}$ requires only one term in the {\it GS}-expansion, which is a huge advantage of the scheme, permitting the use of very few basis terms. Finally, in the limit $\kappa \rightarrow \infty$, the {\it GS}-scheme turns into  the Hermite-spectral scheme, and this due to our specific construction, in particular see \eqref{rec1} and \eqref{rec2}.
\end{remark}
%%%%%%%%%%%%%%%%%%%%%%%%%%%%%%%%%%%%%%%%%%%%%%
\section{Numerical simulations} \label{SEC5}
%%%%%%%%%%%%%%%%%%%%%%%%%%%%%%%%%%%%%%%%%%%%%%%
In this section we present the numerical simulations obtained with the spectral methods constructed in this work, {\it i.e.} the rational Chebyshev ({\it RC}-scheme)  as well as the Gram-Schmidt spectral method ({\it GS}-scheme), and compare them with a standard finite-difference method ({\it FD}-scheme) to solve the energetic Fokker-Planck equation. We are particularly interested in precision, simulation times as well as the asymptotic convergence rates towards the equilibrium, for different values of the parameter $\kappa$. As no analytical solution is available, a reference solution will be constructed via a second-order {\it FD}-scheme on a very fine time and velocity mesh.  All the simulations we performed are carried out in Python $3.11.5$ language, with double precision and using the standard Python libraries, in particular the \textit{mpmath} library for the Whittaker functions \eqref{whitt}.
%%%%%%%%%%%%%%%%%%%%%%%%%%%%%%%%%%%%%%%%%%%	
\subsection{Numerical simulations with a standard FD-method}\label{SEC51}
%%%%%%%%%%%%%%%%%%%%%%%%%%%%%%%%%%%%%%%%%%%	
Let us start with presenting the test case and a finite-difference scheme, aiming to solve numerically the following Fokker-Planck equation
\be \label{FP_num}
(FP)_{\kappa,a}\,\,\,\left\{
  \begin{array}{l}
    \ds  \partial_t f = \partial_v \left[\left(\dfrac{v}{1+\frac{v^{2}}{2\kappa}}+av\right)f+\partial_{v}f\right]\,, \quad \forall (t,v) \in (0,\infty) \times \RR\,,\\[3mm]
    \ds f(0,\cdot)=f_{\mathrm{in}}(\cdot)\,,
  \end{array}
  \right.
\ee
where we shall set either $a=0$ (for the {\it FD}- and {\it RC}-scheme) or $a=10^{-3}$ (regularized version for the {\it GS-}scheme), play with different values of $\kappa \in (1,50)$, and choose an initial condition of two-bump form, composed of two shifted $\kappa$-distribution functions
\be \label{iii}
%f_{\mathrm{in}}(v):={1 \over \sqrt{2\, \pi} \, \sigma}\,e^{-{(v+u)^2 \over 2\, \sigma^2}}+{1 \over \sqrt{2\, \pi} \, \sigma}e^{-{(v-u)^2 \over 2\, \sigma^2}}\,, \qquad \sigma=0.5\,, \quad u=2\,.
f_{\mathrm{in}}(v):={1 \over 2} \,\left[ f_\kappa(v+u)+f_\kappa(v-u)\right]\,,  \quad \forall v \in \RR\,\,\, \textrm{and} \,\,\,  u=2\,.
\ee
In the long-time limit $t \rightarrow \infty$ the solution to this problem tends towards the stationary solution, given by
$$
f_\infty(v):={\mathfrak n}\, f_{\kappa,a}(v)\,, \quad \forall v \in \RR\,; \qquad {\mathfrak n}:={\langle f_{\mathrm{in}} \rangle \over \langle f_{\kappa,a} \rangle}\,,
$$
where we recall the time-decay Theorems \ref{Theorem_Dolbeault} and \ref{Thm_kappa} as well as the notation
$\langle g\rangle:= \int_\RR g(v)\, \dD v$.
Let us remark that the initial condition enters the asymptotic limit function $f_\infty$ only via the mass-constant ${\mathfrak n}$, which has thus to be computed precisely.\\

\noindent To investigate the performance of both spectral schemes, let us introduce a finite-difference scheme as reference method. For this, set $\kappa=3$ in this section and denote for simplicity $\eta(v):= (1+\frac{v^{2}}{2\kappa})^{-1}+a$.
The {\it FD}-scheme is based now on the truncation of the velocity domain to $[-v_{\mathrm{max}},v_{\mathrm{max}}]$,  with $v_{\mathrm{max}}$ to be defined later on, fixing homogeneous Dirichlet boundary conditions for $f$ and introducing a velocity mesh $\{v_{j}\}_{ j=0}^{N_v}$ (not necessarily uniform). Let us also define the mid-points $v_{j+1/2} := (v_j+v_{j+1})/2$ on which equation \eqref{FP_num} will be approximated via a finite-volume scheme. Integrating \eqref{FP_num} over the cell $[v_{j-1/2}, v_{j+1/2}]$ and using the finite-volume notation $f_j(t) := \frac{1}{\Delta v_j} \int_{v_{j-1/2}}^{v_{j+1/2}} f(t, v)\, \dD v$ with cell length $\Delta v_j:=v_{j+1/2}-v_{j-1/2}$, yields
\begin{eqnarray}
\frac{d}{dt} f_j(t) &=&  \frac{1}{\Delta v_j} \Big[  \Big( \eta_{j+1/2}\,v_{j+1/2}\,f(t, v_{j+1/2}) + (\partial_v f)(t, v_{j+1/2})   \Big) \nonumber\\
&& \hspace{3cm}-   \Big( \eta_{j-1/2}\,v_{j-1/2}\,f(t, v_{j-1/2}) + (\partial_v f)(t, v_{j-1/2})   \Big) \Big] \nonumber\\
&\approx &  \frac{1}{\Delta v_j}  \Big[  \Big( \eta_{j+1/2}\,v_{j+1/2}\,\frac{f_j(t) + f_{j+1}(t)}{2} + \frac{f_{j+1}(t) - f_{j}(t)}{v_{j+1} - v_j}  \Big) \nonumber\\
&& \hspace{3cm}-   \Big( \eta_{j-1/2}\,v_{j-1/2}\frac{f_j(t) + f_{j-1}(t)}{2} +  \frac{f_j(t) - f_{j-1}(t)}{v_{j} - v_{j-1}}   \Big) \Big] \,,
%\nonumber\\
%&=& \frac{1}{\Delta v_j} \Big\{  f_{j-1}(t)  \Big( - \eta_{j-1/2}\,{v_{j-1/2} \over 2} + \frac{1}{v_{j} - v_{j-1}}    \Big)  \nonumber\\ 
%&& +  f_{j}(t)   \Big[  \left( \eta_{j+1/2}\,{v_{j+1/2}\over 2}- \frac{1}{v_{j+1} - v_{j}}\right) - \left(\eta_{j-1/2}\, {v_{j-1/2} \over 2} + \frac{1}{v_{j} - v_{j-1}} \right) \Big] \nonumber\\ 
%\label{discrete_Q}
%&& + f_{j+1}(t)  \Big(  \eta_{j+1/2}\,\frac{v_{j+1/2}}{2} + \frac{1}{v_{j+1} - v_{j}}    \Big) \Big\}\,, 
\end{eqnarray}
for all $j=1, \cdots, N_v-1$, whereas we set for the boundary values $f_0=f_{N_v} \equiv 0$.\\

In what concerns the time discretization, an implicit Euler or Crank-Nicolson scheme is used to solve the thus obtained ODE system, namely
$$
({\mathbb I}+\Delta t\, A)\,f^{n+1}=f^n \quad \textrm{or} \quad ({\mathbb I}+{\Delta t \over 2}\, A)\,f^{n+1}=({\mathbb I}-{\Delta t \over 2}\,A)\,f^n\,, \qquad \forall n \in \NN\,,
$$
with $A$ a tridiagonal matrix of size $(N_v-1)\times (N_v-1)$, and starting point $f_j^0:=f_{\mathrm{in}}(x_j)$.\\

\begin{figure}[ht]
\begin{center}
  \includegraphics[scale=0.4]{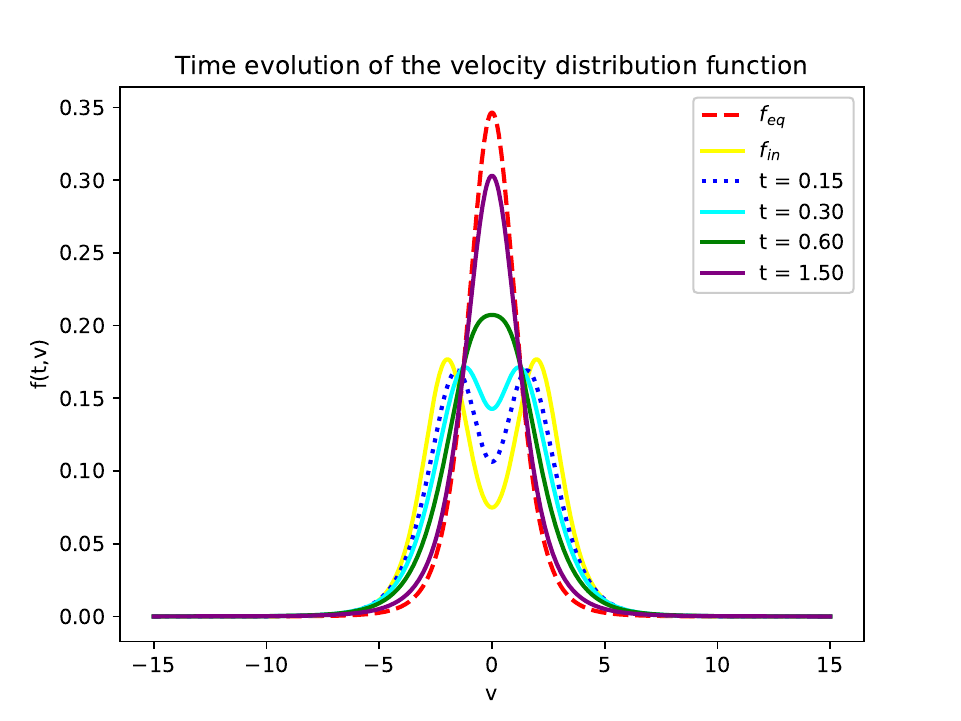}\hfill
  \includegraphics[scale=0.4]{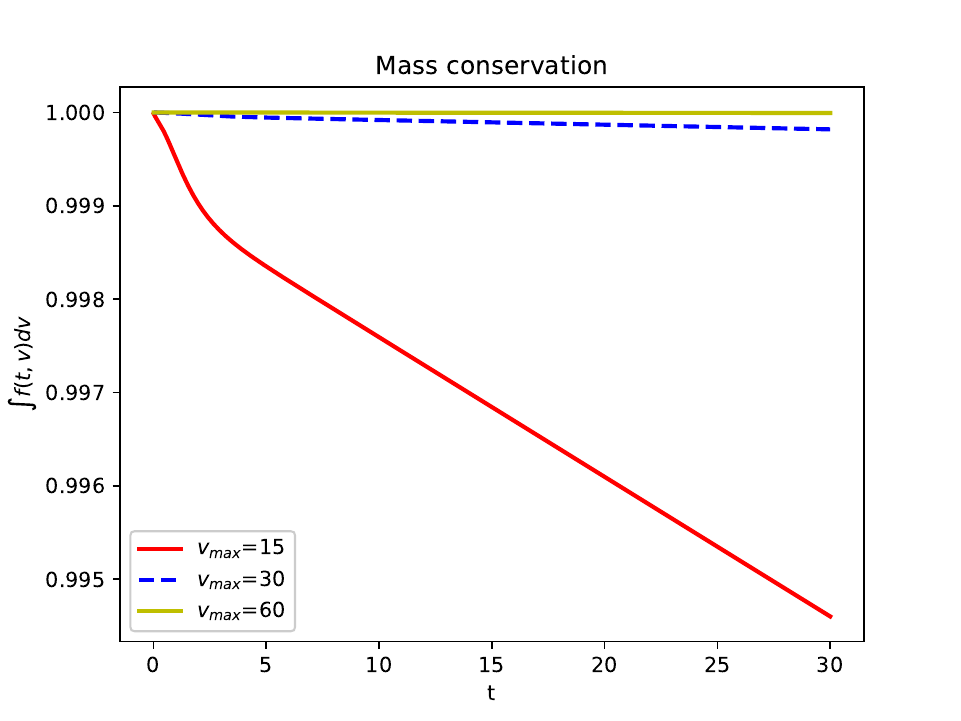}
\end{center}
\caption{{\it FD}-scheme. Left: Time evolution of $f(t,v)$, starting from the two-bump initial condition \eqref{iii} and $\kappa=3$. Right: Plot of the mass $\int_\RR f(t,v)\, \dD v$ for three different truncations $v_{\mathrm{max}}$.}\label{INIT}
\end{figure}

\begin{figure}[ht]
\begin{center}
  \includegraphics[height=4cm, width=0.48\textwidth]{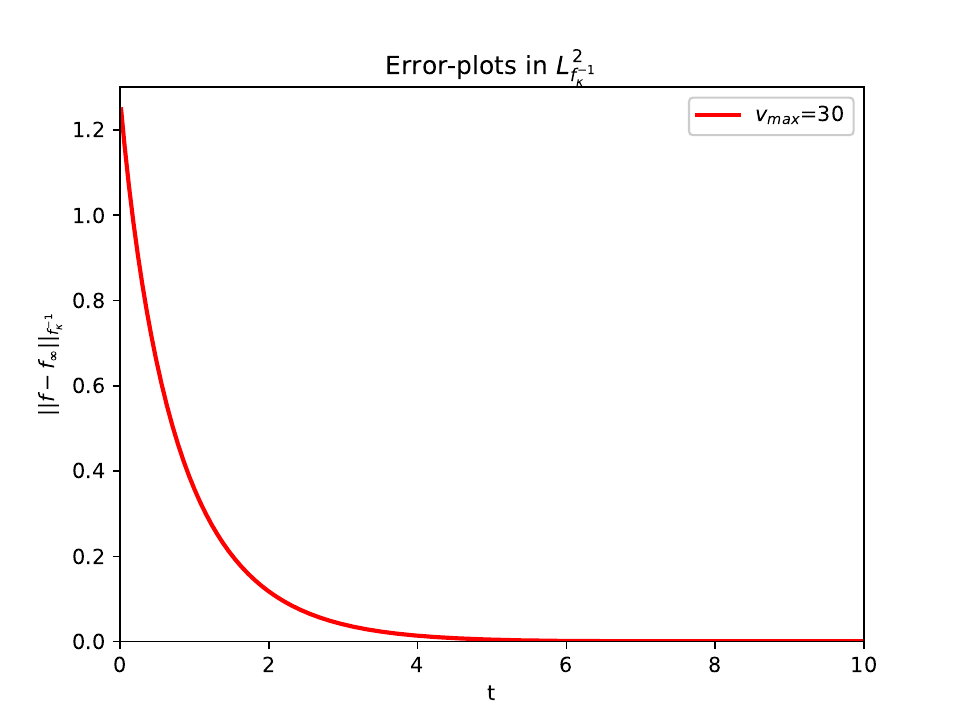}\hfill
  \includegraphics[height=4cm, width=0.48\textwidth]{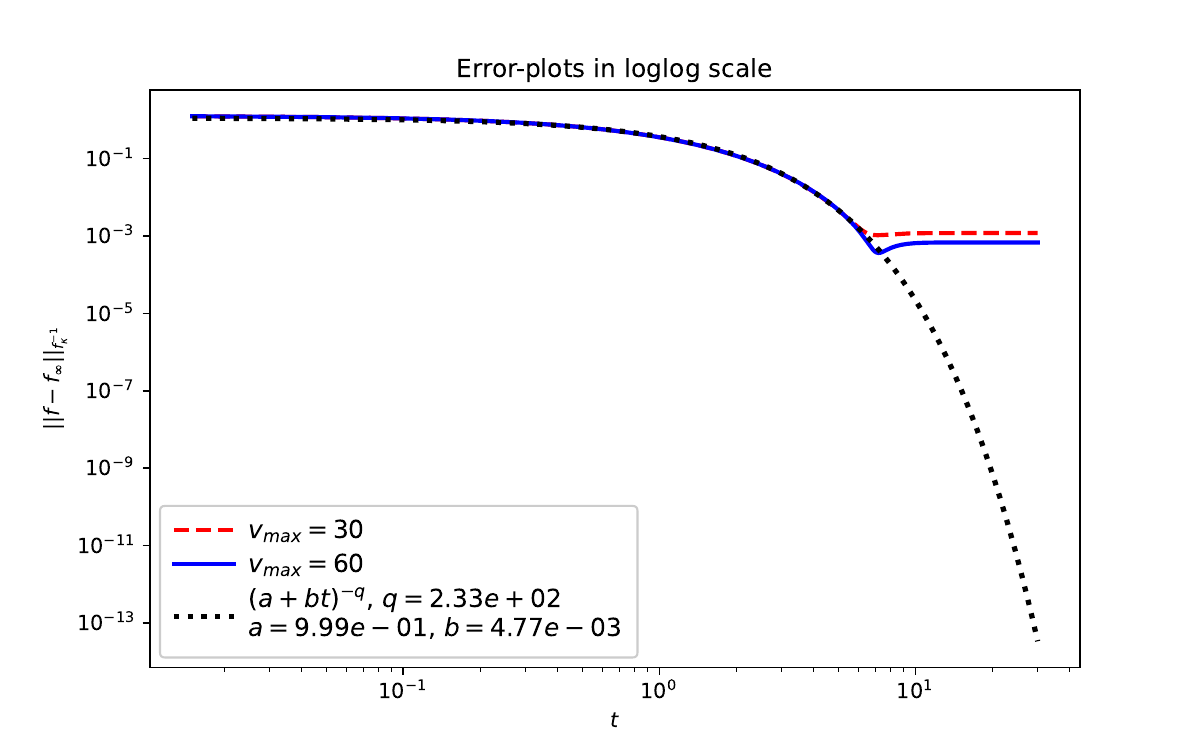}
\end{center}
\caption{{\it FD}-scheme. Time evolution of $||f(t)-f_{\infty}||_{L^2_{f^{-1}_\kappa}}$ for $\kappa=3$. Left: $v_{\mathrm{max}}=30$. Right: Error plots in log-log scale and for two different domain truncations $v_{\mathrm{max}}=30$ and $v_{\mathrm{max}}=60$, as well as fitting curve $g(t)=(a+b\, t)^{-q}.$}\label{Error}
\end{figure}

%, given by
%\begin{eqnarray*}
%A_{j,j-1}&:=&  -\frac{1}{\Delta v_j} \Big( - \eta_{j-1/2}\, \frac{v_{j-1/2}}{2} + \frac{1}{v_{j} - v_{j-1}} \Big)\,,\nonumber\\
%A_{j,j}&:=&  -\frac{1}{\Delta v_j} \Big[ \left( \eta_{j+1/2}\,{v_{j+1/2}\over 2}- \frac{1}{v_{j+1} - v_{j}}\right)  - \left(\eta_{j-1/2}\, {v_{j-1/2} \over 2} + \frac{1}{v_{j} - v_{j-1}} \right)\Big]\,,\nonumber\\
%A_{j,j+1} &:=& - \frac{1}{\Delta v_j}  \Big(  \eta_{j+1/2}\,\frac{v_{j+1/2}}{2} + \frac{1}{v_{j+1} - v_{j}}    \Big). 
%\end{eqnarray*}
%A second-order Crank-Nicolson scheme could be used to improve the accuracy, which is however not the main objective for our comparison studies.\\
To illustrate the evolution of the velocity distribution function, we used the first-order {\it FD}-scheme and plotted on the left of Fig. \ref{INIT} the initial condition $f_{\mathrm{in}}(v)$ corresponding to the two-bump distribution function \eqref{iii}, as well as its evolution towards the asymptotic state $f_\infty(v)$. In order to evaluate the convergence rate towards this equilibrium, we plotted on Fig. \ref{Error} the distance $||f(t)-f_\infty||_{L^2_{f^{-1}_\kappa}}$ as a function of time  (left plot), as well as its $log-log$ scale version (right plot). One remarks that after an initial rapid decay, the error saturates, as observed from the right plot of Fig. \ref{Error},  effect due, among other things, to the truncation of the velocity domain. Moreover, in order to evaluate the algebraic decay rate in time, we observe that the $log-log$ scale curve fits very well with a function of the type $g(t)=(a+b\, t)^{-q}$, where we found $a\approx 1$, $b \approx 5\cdot 10^{-3}$ and $q \approx 233$. Examining the error estimate in 
Theorem \ref{Theorem_Dolbeault}, we must have $-p/4=-q$, meaning thus that $p \approx 930$, which indicates that the scheme is performing much better than it should, after Theorem \ref{Theorem_Dolbeault}. However, let us underline here that it is very intricate to verify numerically the theoretical decay rate, due to several reasons. Firstly, the result of this theorem is not necessarily optimal, the explicit expression of the constant $K_{p,\kappa}$ is not known. Secondly, our numerical simulations are performed with a specific initial condition. And finally, an essential argument is that we truncated the domain in our numerical simulations, a fact that changes an algebraic decay into an exponential decay.  This can be easily seen from the values of $b \ll 1$ and $q \gg 1$,  recalling that $e^{-c\, t}= \lim_{q \rightarrow \infty} \left(1+ {c\, t \over q}\right)^{-q}$.

Finally, let us also observe that the mass is not perfectly conserved numerically (see right of Fig. \ref{INIT}). The reason for this comes from the fact, that we truncated the velocity domain at $v_{\mathrm{max}}=15$ and set homogeneous Dirichlet boundary conditions $f(t,\pm v_{\mathrm{max}})=0$ for all times $t \ge 0$. However, since we are dealing with energetic particles, {\it i.e.} distribution functions with heavy tails, truncating the domain necessarily leads to mass losses. Truncating the domain more far away reduces this error, however increases naturally the computational time. It is because of such difficulties in the numerical discretization (via {\it FD}-schemes) of the energetic particle dynamics, that we proposed in this paper the two spectral schemes, the {\it RC}- and the {\it GS}-scheme. 

%%%%%%%%%%%%%%%%%%%%%%%%%%%%%%%%%%%%%%%%%%%	
\subsection{Numerical results obtained with our spectral methods}\label{SEC52}
%%%%%%%%%%%%%%%%%%%%%%%%%%%%%%%%%%%%%%%%%%%
Let us investigate now the numerical results obtained for problem \eqref{FP_num}-\eqref{iii} with both spectral methods, the {\it RC}-scheme presented in Section \ref{SEC3} and the {\it GS}-scheme designed in Section \ref{SEC4}. The solution $f(t,v)$ of this problem is approximated by one of the two following truncated series (corresponding to the two different basis-sets previously introduced)
$$
f_{RC}(t,v):=\sum_{k=0}^{{\tilde N_{RC}}} {\mathfrak a}_{2k}(t)\, \Theta_{2k}(v)\,, \qquad 
f_{GS}(t,v):=\sum_{k=0}^{N_{GS}} \alpha_k(t)\, p_k(v)\, f_{\kappa,a}(v)\,,
$$
where ${\tilde N_{RC}}=N_{RC}/2$ and $N_{RC},\,N_{GS}$ are given in Table \ref{tab1}, whereas the spectral coefficients ${\mathfrak a}_{2k}(t)$ (resp. $\alpha_k(t)$) are solutions of the linear systems \eqref{L_syst} (resp. \eqref{ODE}) for the {\it RC}-scheme (resp. {\it GS}-scheme), solved via a first-order implicit Euler-method. We omitted for simplicity reasons to index the functions $f_{RC}$ and $f_{GS}$ by the truncation index. Let us recall that in the asymptotic limit $t \rightarrow \infty$ the steady state needs for the {\it RC}-scheme ${\kappa+1 \over 2}$ terms, thus ${\tilde N_{RC}}={\kappa-1 \over 2}$, whereas the {\it GS}-scheme requires only one term, thus $N_{\textrm{GS}}=0$. The initial condition we choose is again the two-bump distribution given in \eqref{iii}, its time-evolution being similar to the one plotted on the left of Fig. \ref{INIT}.

\begin{table}[htb] 
\caption{Numerical parameters used for the confrontation of the different schemes.}
\begin{tabular}{|c||c||c|c|c|}
\hline 
    & ${\it FD}_{\mathrm{ref}}$-scheme  & {\it FD}-scheme & {\it RC}-scheme  & {\it GS}-scheme \\
\hline
 \hline
 Discret. &Crank-Nicolson&	impl. Euler & impl. Euler& impl. Euler \\
 \hline
 $a$ &0&	0 & 0& $10^{-3}$ \\
 \hline
 $\kappa$ &$3$ or $31$&	$3$ or $31$ & $3$  & $31$   \\
 \hline
 $N_v$; $v_{\mathrm{max}}$ &$10\,001$; $25$& $[101\cdots 501]$; $15$& $\minuso$ &$\minuso$  \\
 \hline
 $N_{RC/GS}$ & $\minuso$ & $\minuso$ & $8$, $10$, $12$, $16$ & $6$, $8$, $10$, $12$, $16$ \\
 \hline
 $T$; $N_t$ & $2$;\, $2\cdot10^3$ & $2$\,; $[50 \cdots 10^3]$ &$2$; $[50 \cdots 10^3]$  &$2$\,; $[50 \cdots 10^3]$  \\
\hline

\end{tabular}
 \label{tab1}
\end{table}

\begin{figure}[ht]
\begin{center}
  \includegraphics[scale=0.4]{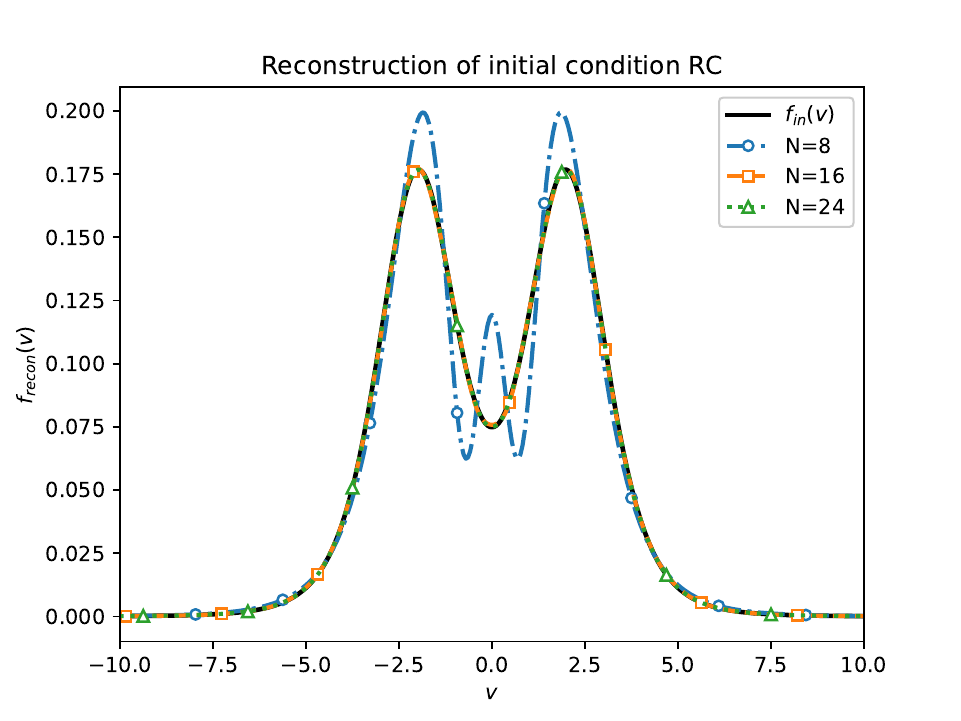}\hfill
  \includegraphics[scale=0.4]{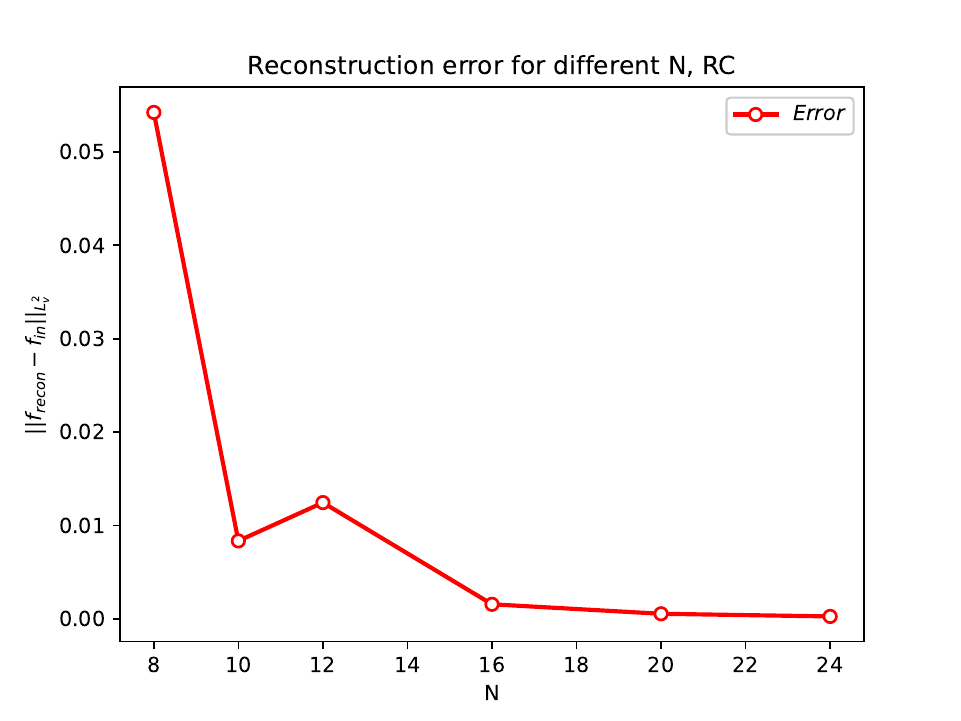}
\end{center}
\caption{Reconstruction of the initial condition $f_{\mathrm{in}}(v)$ via the $RC$ basis-set ($\kappa=3$), for several truncation indices $N_{RC}$. Left: Velocity distribution function. Right: Error plots $||f_{recon}-f_{\mathrm{in}}||_{L^2_v}$. For $N=16$, the error is of $1.6\cdot 10^{-3}$.}\label{reconstr_RC}
\end{figure}

\begin{figure}[ht]
\begin{center}
  \includegraphics[scale=0.4]{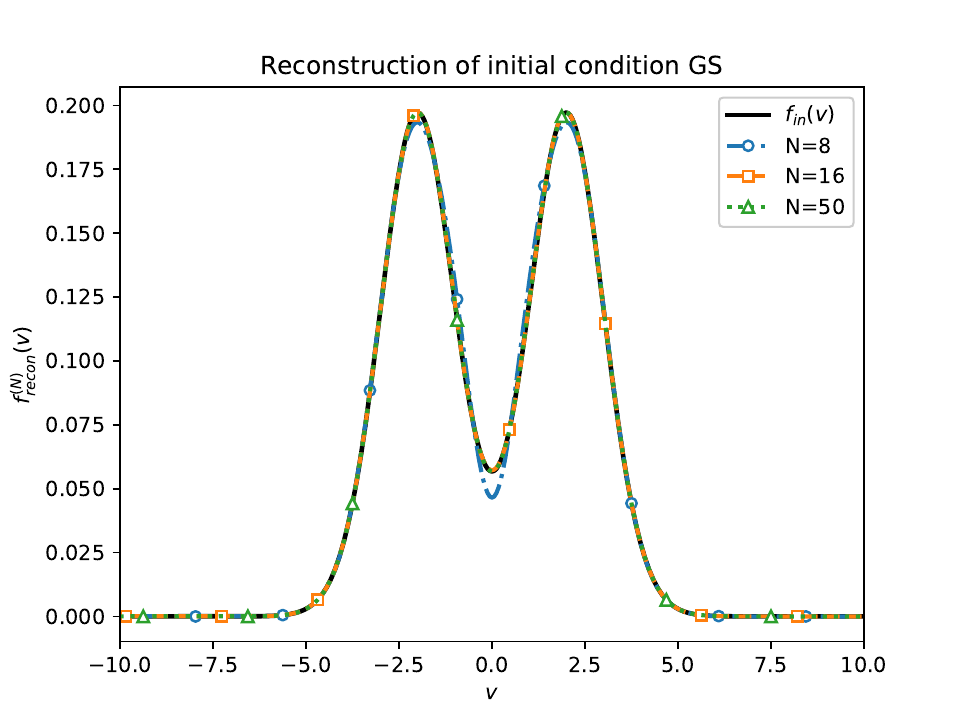}\hfill
  \includegraphics[scale=0.4]{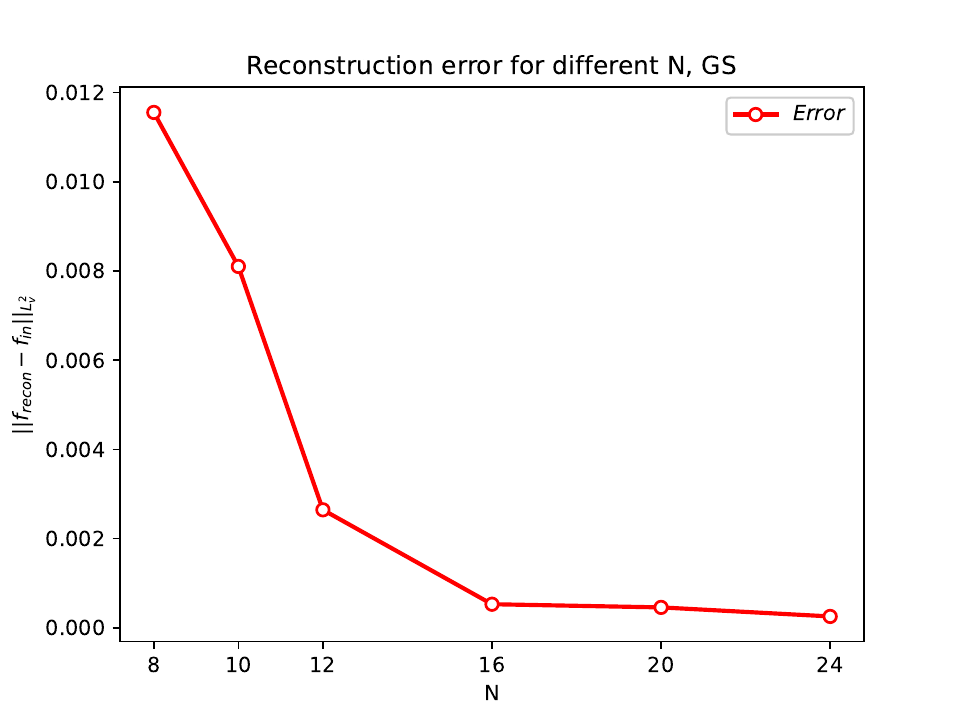}
\end{center}
\caption{Reconstruction of the initial condition $f_{\mathrm{in}}(v)$ via the $GS$ basis set ($\kappa=31$), for several truncation indices $N_{GS}$. Left: Velocity distribution function. Right: Error plots $||f_{recon}-f_{\mathrm{in}}||_{L^2_v}$.  For $N=12$, the error is of $2.6\cdot 10^{-3}$.} \label{reconstr_GS}
\end{figure}

To start the procedure, we need first to project the initial condition on the corresponding orthogonal basis set, which is done by computing
$$
{\mathfrak a}_{2\ell}(0)={ 1 \over {\mathfrak c}_\ell}\, \int_\RR f_{\mathrm{in}}(v)\, \Theta_{2\ell}(v)\, f_\kappa^{-1}(v)\, \dD v  \,, \quad \alpha_\ell(0)={1 \over \gamma_\ell}\, \int_\RR f_{\mathrm{in}}(v)\, p_\ell(v)\, \dD v  \,, \quad \forall \ell \ge 0\,,
$$
where ${\mathfrak c}_0={\pi}$, ${\mathfrak c}_\ell={\pi/2}$ for $\ell \ge 1$ and $\gamma_\ell$ is the weighted norm $\gamma_\ell=||p_\ell ||^2_{f_{\kappa,a}}$. The computation of these integrals is an important step in our spectral methods, a bad (not-accurate or too time-consuming) computation would lead to a loss of all the advantages of the spectral methods. Standard integration methods are used here, and we illustrated, for validation, in Fig. \ref{reconstr_RC}  and \ref{reconstr_GS} the reconstruction of the initial condition in the corresponding basis set. From these  figures one can identify the number of basis-functions required in order to have a precise simulation in the initial time-phase, for example the {\it RC}-scheme requires $N_{RC}=16$ for a precision of $1.6\cdot 10^{-3}$ whereas the {\it GS}-scheme requires $N_{GS}=12$ basis functions for a precision of $2.6\cdot 10^{-3}$.\\

A few simulations are carried out now in the aim to compare the proposed schemes with respect to convergence, accuracy and computational efficiency. For these comparison studies, a reference solution $f_{\textrm{ref}}(t,v)$ of \eqref{FP_num}-\eqref{iii} is constructed via a second-order {\it FD}-scheme (Crank-Nicolson time-discretization) on a very large domain and with a very fine mesh (see Table \ref{tab1}). This reference solution can be considered somehow as an exact solution. The discrete $L^2_v$-norm we shall use in the following for the error investigations, is defined as
$$
||f_{\mathrm{ref}}-f_{\mathrm{num}}||^2_{L^2_v}(t_\star):= { 2\, v_{\mathrm{max}} \over N_v}\sum_{j=0}^{N_v-1} |f_{\mathrm{ref}}(t_\star,v_j)-f_{\mathrm{num}}(t_\star,v_j)|^2\,.
$$

%Depending on the scheme, we shall speak about $N$-convergence, when convergence is achieved for the spectral method, by increasing the degree of the polynomials or equally the expansion truncation index $N$, whereas we shall speak about $\Delta v$-convergence, when convergence is achieved for the FD-scheme by increasing the number of grid points $N_v$.\\

\begin{figure}[ht]
\begin{center}
  \includegraphics[scale=0.4]{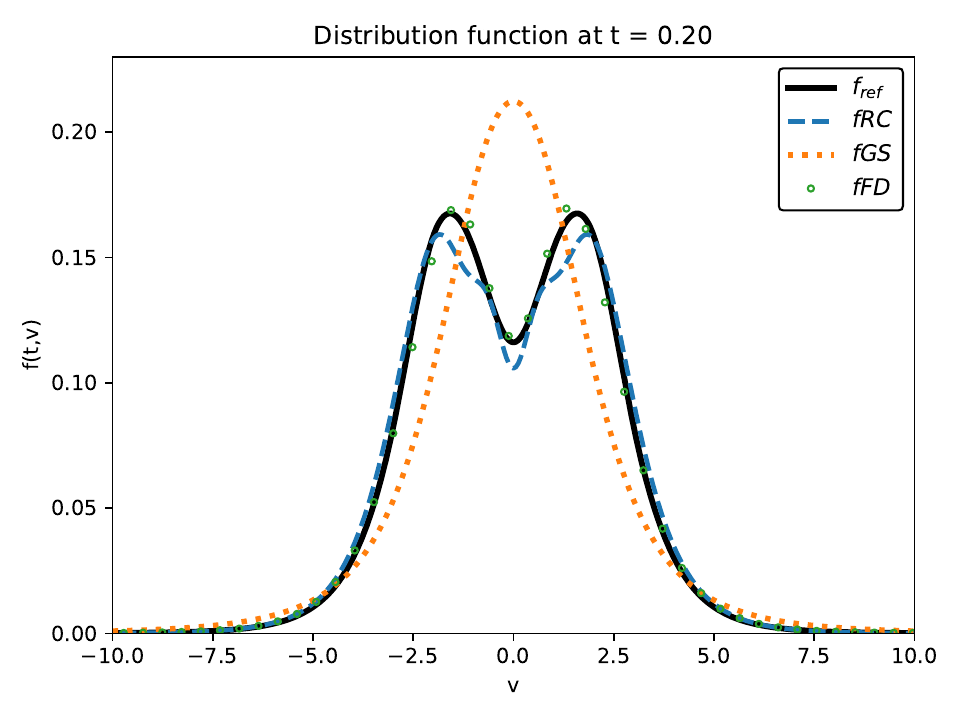}\hfill
  \includegraphics[scale=0.4]{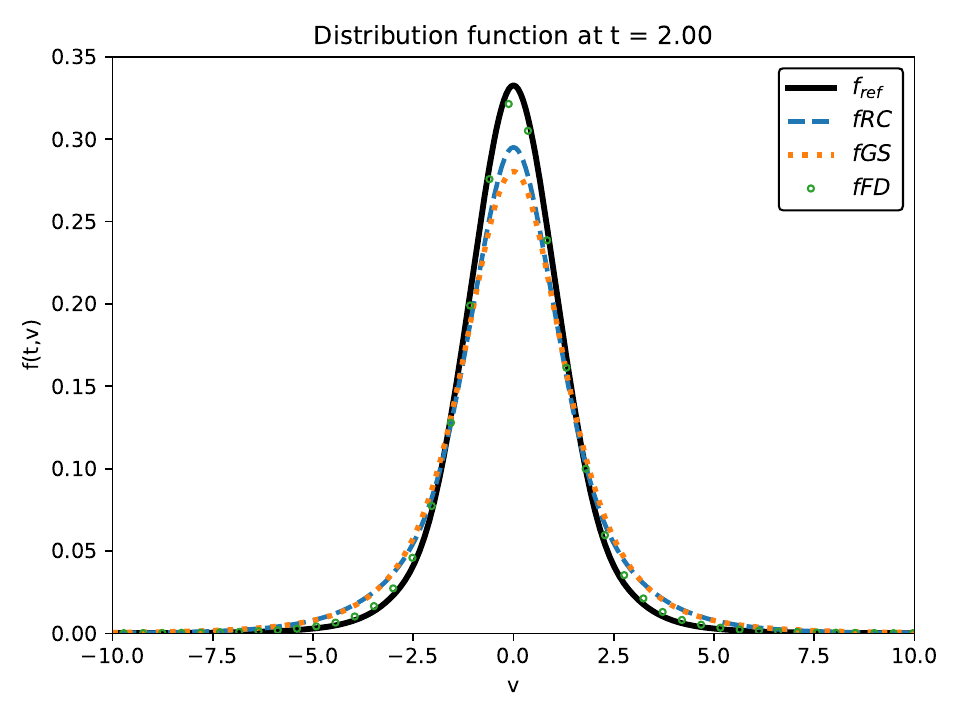}
\end{center}
\caption{Left: Velocity distribution functions $f_{\mathrm{ref}}(t_\star,v)$  and $f_{\mathrm{num}}(t_\star,v)$ for $t_\star=0.2$ and $\kappa=3$. Right: Same velocity distribution functions at instant $t_\star=2$. Furthermore $N_{RC}=N_{GS}=10$. }\label{f_k3}
\end{figure}

\begin{figure}[ht]
\begin{center}
  \includegraphics[scale=0.4]{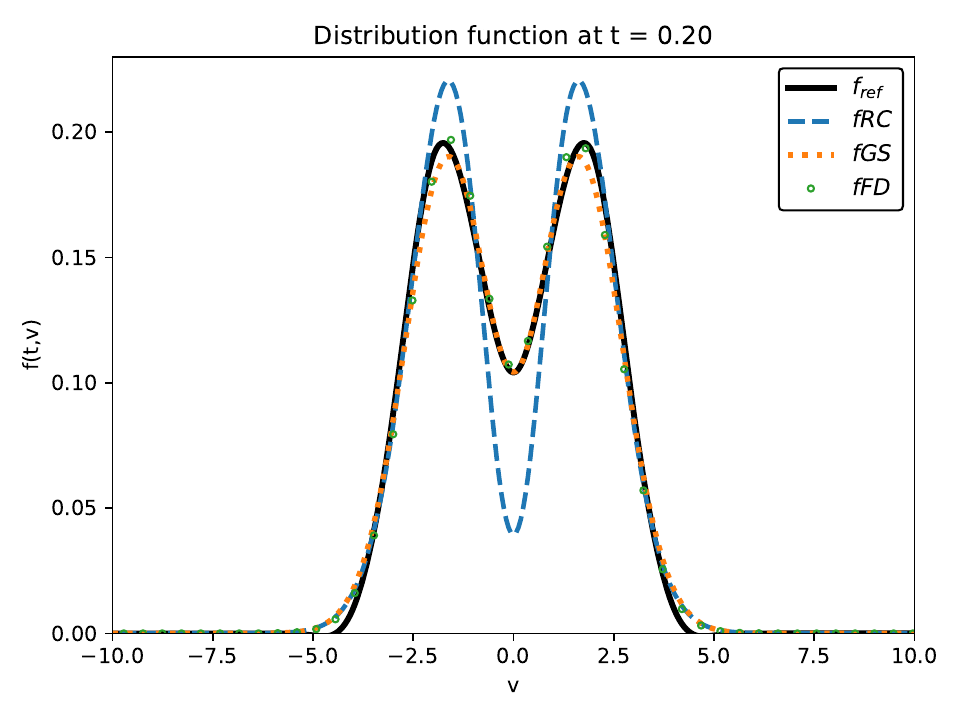}\hfill
  \includegraphics[scale=0.4]{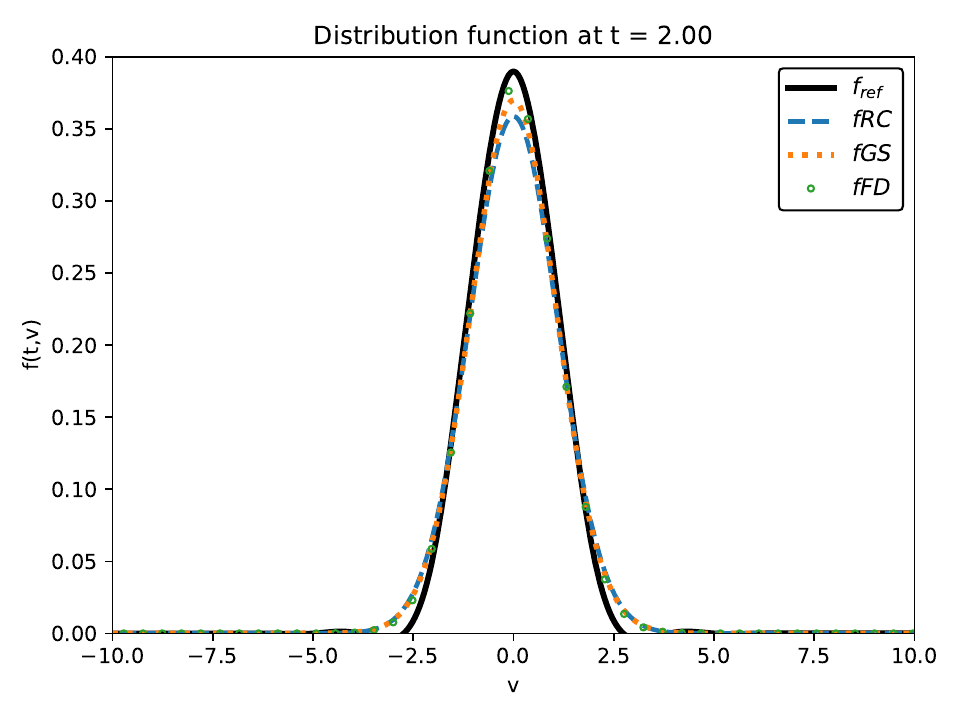}
\end{center}
\caption{Left: Velocity distribution functions $f_{\mathrm{ref}}(t_\star,v)$ and $f_{\mathrm{num}}(t_\star,v)$ for $t_\star=0.2$ and $\kappa=31$. Right: Same velocity distribution functions at instant $t_\star=2$. Again $N_{RC}=N_{GS}=10$.}\label{f_k31}
\end{figure}

Let us start by presenting in Fig. \ref{f_k3} and \ref{f_k31} the curves corresponding to the reference distribution function $f_{\mathrm{ref}}(t,v)$ as well as the three distributions $f_{\mathrm{num}}(t,v)$ computed with the first-order ${\it FD}$-, ${\it RC}$- and ${\it GS}$-schemes at two different times $t_\star=0.2$ and $t_\star=2$ as well as for two different $\kappa$-values, $\kappa=3$ and $\kappa=31$. What can be observed from these plots is that for small times, the accurate projection of the initial condition on the corresponding basis-set is of great importance. In particular the $RC$-scheme has some problems for large $\kappa$ and this due to the necessity to take a large number of basis-functions for the reconstruction, namely $N_{RC} > \kappa-1$ as explained in Remark \ref{remCheb}, whereas the $GS$-scheme has some problems for small $\kappa$-values, due to the hot tails of the distribution functions and the need to take a large (non-physical) regularization parameter $a>0$ in order to overcome stability problems, as explained in Remark \ref{remGS}.
For larger times, all numerical solutions get closer and closer to the equilibrium solution $f_\infty(v)$, and the initial condition no longer plays an important role (except in mass-computation). Thus, each scheme is better adapted for a particular $\kappa$-regime. We shall therefore take in the following $\kappa=3$ for the $RC$-scheme and $\kappa=31$ for the $GS$-scheme. As explained earlier, one can read from Fig. \ref{reconstr_RC}  and \ref{reconstr_GS} the number of basis-functions required for each scheme to have a good precision, especially in the small-time regime.\\



\color{black}

\begin{figure}[ht]
\begin{center}
  \includegraphics[scale=0.4]{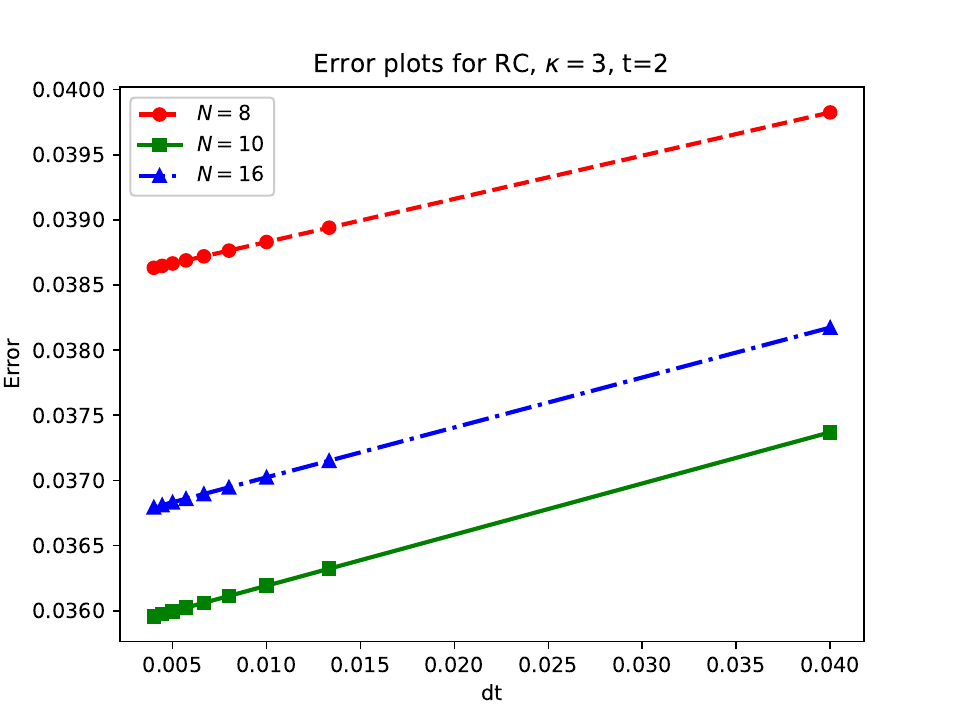}\hfill
  \includegraphics[scale=0.4]{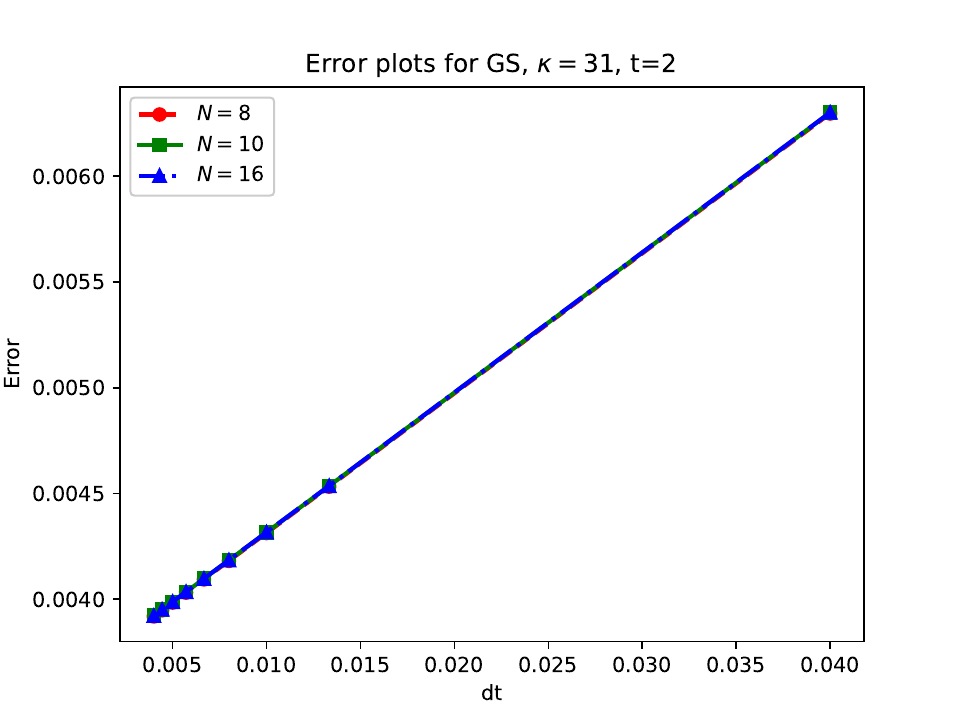}
\end{center}
\caption{Error curves $||f_{\mathrm{ref}}(t_\star,\cdot)-f_{\mathrm{num}}(t_\star,\cdot)||_{L^2_v}$ for $t_\star=2$ with {\it RC}-scheme (left, $\kappa=3$) and {\it GS}-scheme (right, $\kappa=31$). The reference sol. ({\it FD}-scheme) is discretized via Crank-Nicolson.}
\label{conv_RC}
\end{figure}

To demonstrate the convergence of our spectral schemes, we depicted on Fig. \ref{conv_RC} the $L^2_v$-errors of our spectral methods with respect to the reference solution, {\it i.e.} $||f_{\mathrm{ref}}(t_\star,\cdot)-f_{RC/GS}(t_\star,\cdot)||_{L^2_v}$, errors as a function of the time-discretization step $\Delta t$, at time instant $t_\star=2$ and for three different truncation indices $N_{RC/GS}=8,10,16$. For the Chebyshev scheme we used $\kappa=3$, while for the Gram-Schmidt scheme we chose $\kappa=31$ in accordance with the remarks made above. What can be observed is that the errors are indeed small, which signifies that both spectral schemes approximate well the unknown solution. As Euler's implicit scheme is used to advance in time in the spectral methods, we expect a linear convergence in time, which is indeed visible from these two plots, however the slopes are very small, of order $10^{-2}$, meaning the convergence is very slow, or even saturates somehow. This discordance is simply due to the fact that the reference solution is not an exact solution, but a numerical solution, with fixed time-step of $\Delta t=10^{-3}$ and homogeneous Dirichlet boundary conditions, thus containing also numerical errors. Our ${\it RC}$- and ${\it GS}$-schemes seem to become rapidly too accurate, to be anymore compared with this "reference solution". A higher-order ($4^{th}$-order) {\it FD}-scheme as well as transparent boundary conditions could be a good idea, for a better comparison and should be investigated in future works. We also observe that the three lines, corresponding to the truncation indices, $N_{RC}$ are not in the expected order. This comes from the fact that with more polynomials in the expansion, more computations are carried out, which can indeed lead to an accumulation of the round-off errors. Thus an optimum has to be found between small truncation errors and small accumulation errors.
%the convergence of both schemes, the errors being linear in the discretization step $\Delta t$, which is completely natural, as we are discretizing the system in time via the implicit Euler scheme. Note that the slope becomes very small in the last plot of both Fig. \ref{conv_RC} resp. \ref{conv_GS}. This is simply due to saturation effects, the reference solution being not an exact solution, but a numerical solution, where we fixed the time-step to $\Delta t=10^{-3}$. This artifact is clearly visible in the last plot of Fig. \ref{conv_GS}, where the error saturates at some value of $1.212 \cdot 10^{-3}$. As a consequence, our ${\it RC}$- and ${\it GS}$-schemes become at large times too accurate, to be anymore compared with this reference solution.
In the ${\it GS}$-scheme the number of basis-functions $N_{GS}$ does not play anymore a role, the three lines become superposed after some transition-phase in time, the asymptotic regime requiring only one term in the expansion.\\

\begin{figure}[ht]
\begin{center}
  \includegraphics[height=4cm,width=0.48\textwidth]{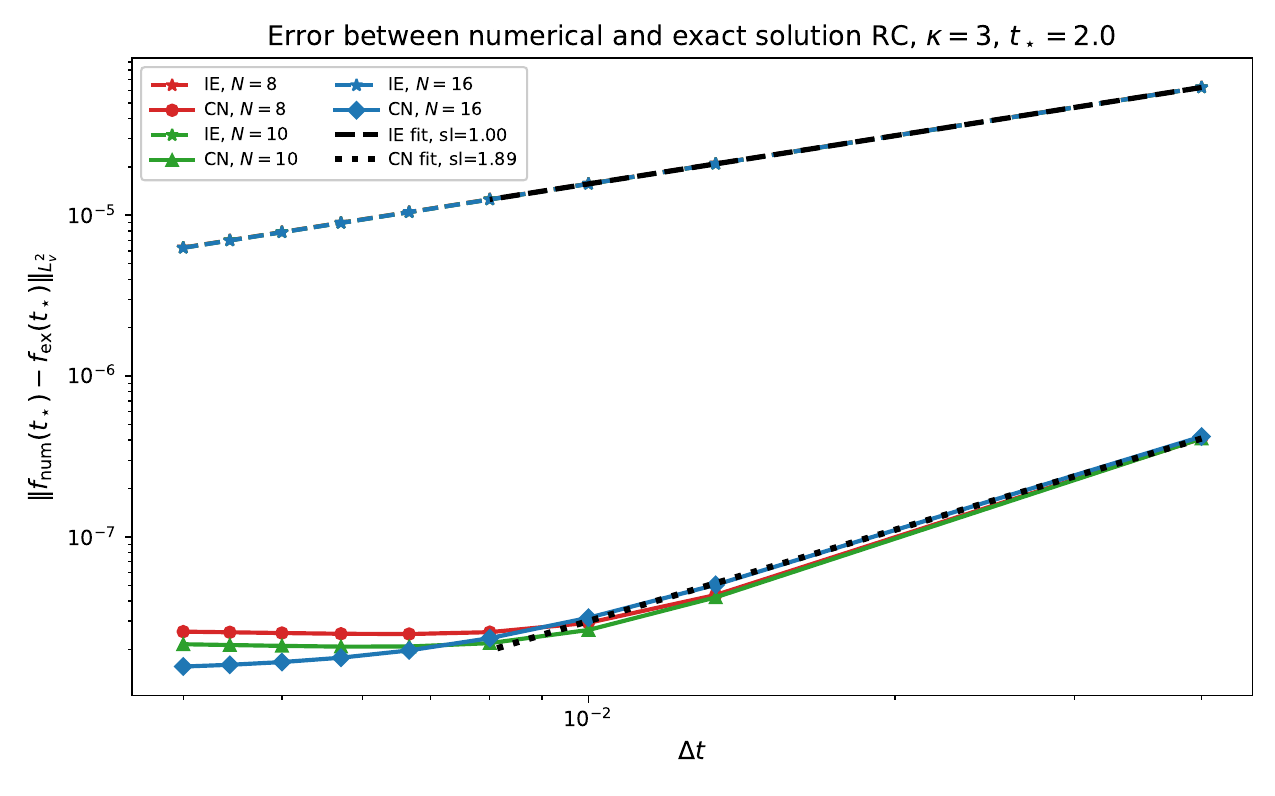}\hfill
  \includegraphics[height=4cm,width=0.48\textwidth]{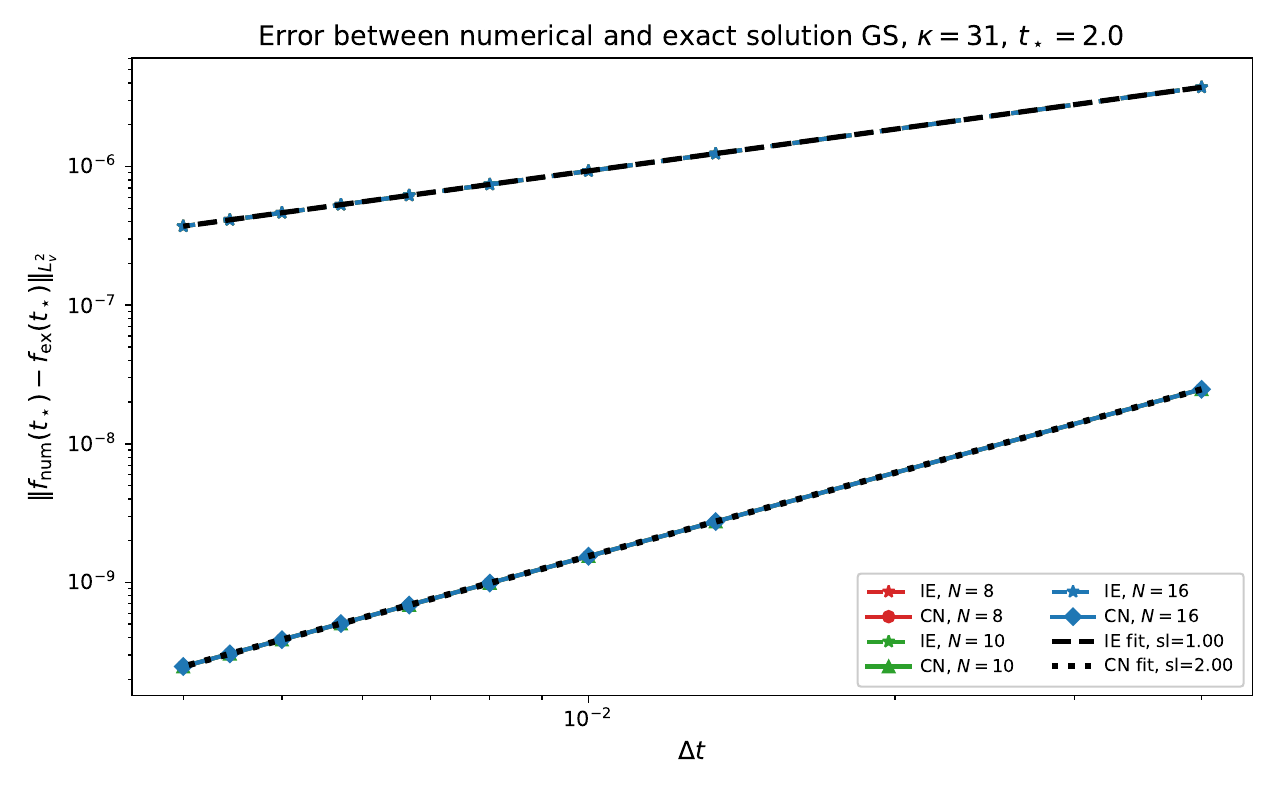}
\end{center}
\caption{Error curves $||f_{\mathrm{ex}}(t_\star,\cdot)-f_{\mathrm{num}}(t_\star,\cdot)||_{L^2_v}$ for $t_\star=2$ with {\it RC}-scheme (left, $\kappa=3$) and {\it GS}-scheme (right, $\kappa=31$). Here $f_{\mathrm{ex}}$ denotes the manufactured solution \eqref{start} of the forced Fokker-Planck eq. \eqref{FP_num}, and $f_{\mathrm{num}}$ the spectral sol. discretized either via implicit Euler or Crank-Nicolson.
%Relative error curves $||f_{dt}(t_\star,\cdot)-f_{dt/2}(t_\star,\cdot)||_{L^2_{f^{-1}_\kappa}}/||f_{dt/2}(t_\star,\cdot)||_{L^2_{f^{-1}_\kappa}}$ wrt. $dt$, in log-log-scale, for $t_\star=2$ and $N_{RC/GS}=10$. Left: RC-scheme, $\kappa=3$; Right: GS-scheme, $\kappa=31$.
}\label{auto_conv}
\end{figure}

To overcome the fact that no analytical solution is available for the verification of the convergence of our spectral schemes, we decided to construct a manufactured (exact) solution, which is the exact solution of a slightly modified Fokker-Planck equation, containing an additional forcing term.  This constructed manufactured solution (even if physically not relevant) can now be used as a benchmark solution for verification purposes. Both spectral methods are now used to solve this modified Fokker-Planck equation, and we compare in Fig. \ref{auto_conv} the obtained numerical solutions (via a first-order implicit Euler as well as a  second-order Crank-Nicolson time-discretization), with the exact manufactured solution, fact which permits finally to show the convergence of our spectral methods, in particular we find the expected slopes for the Euler as well as the Crank-Nicolson method.\\

To finish, let us present here the construction of the manufactured solution corresponding to the comparison in Fig. \ref{auto_conv}. We started with
\be \label{start}
f_{\rm ex}(t,v)=f_{\mathrm{eq}}(v)\left(1+\varepsilon e^{-t}g(v)\right), \qquad \forall v \in \RR\,,
\ee
where $\eps =0.1$ and 
\[
g(v)=C_2^L(v)-\overline C_2,
\qquad
C_2^L(v)=\frac{v^2-2\kappa}{v^2+2\kappa},
\qquad 
\overline C_2
=
\frac{\int_{\mathbb R} f_{\mathrm{eq}}(v)C_2^L(v)\,\dd v}
{\int_{\mathbb R} f_{\mathrm{eq}}(v)\,\dd v}.
\]
Inserting this specific function \eqref{start} in our Fokker-Planck equation \eqref{FP_num}, permits to obtain analytically a forcing (source) term, namely $f_{\rm ex}$ is the exact solution of the following forced Fokker-Planck equation 
\be \label{forced}
\partial_t f
=
\partial_v\left(
f_{\mathrm{eq}}(v)\partial_v\left(\frac{f}{f_{\mathrm{eq}}}\right)\right)
+
e^{-t}\, {\mathcal S}(v), \qquad \forall (t,v) \in (0,\infty) \times \RR\,,
\ee
with
\[
{\mathcal S}(v)
=
-\varepsilon f_{\mathrm{eq}}(v)
\left[
g(v)+g''(v)+(\log f_{\mathrm{eq}})'(v)g'(v)
\right].
\]
For the $RC$-scheme, one has
\[
f_{\mathrm{eq}}(v)=f_\kappa(v)
=
c_\kappa\left(1+\frac{v^2}{2\kappa}\right)^{-\kappa},
\qquad
(\log f_{\mathrm{eq}})'(v)
=
-\frac{v}{1+\frac{v^2}{2\kappa}},
\]
whereas for the $GS$-scheme, one has
\[
f_{\mathrm{eq}}(v)=f_{\kappa,a}(v)
=
c_\kappa\left(1+\frac{v^2}{2\kappa}\right)^{-\kappa}e^{-av^2/2},
\qquad
(\log f_{\mathrm{eq}})'(v)
=
-\frac{v}{1+\frac{v^2}{2\kappa}}-av.
\]

%\begin{figure}[ht]
%\begin{center}
%  %\includegraphics[scale=0.3]{Plots/error_GS_k31_t02.eps}\hfill
%  \includegraphics[scale=0.3]{Plots/error_GS_k31_t2.pdf}\hfill
%  \includegraphics[scale=0.3]{Plots/error_GS_k31_t10.eps}
%\end{center}
%\caption{Error curves $||f_{\mathrm{ref}}(t_\star,\cdot)-f_{GS}(t_\star,\cdot)||_{L^2_v}$ for $\kappa=31$, three time instants $t_\star=0.2;\,\,2;\,\,\, 10$ and three different truncation indices $N=8;\,\, 10;\,\, 16$.}\label{conv_GS}
%\end{figure}

Having proven the convergence of both of our spectral schemes, let us come back now again to the comparison of our spectral methods with the standard finite-difference scheme, in particular let us investigate whether the spectral schemes have an advantage with respect to standard finite-difference methods.
Naturally, increasing the truncation domain via $v_{\mathrm{max}}$, the number of grid points $N_v$ or the truncation index $N_{RC/GS}$ shrinks the error, but also increases the simulation time. To compare the computational times of the three different methods, we performed several simulations for different sets of parameters, fixing a large final time $t_\star=10$. In Table \ref{tab_t10} we collect the errors obtained (with respect to the manufactured exact solution of \eqref{forced}) in mass, indicated by ${\EuFrak{E}_m}$, the total error $||f_{\mathrm{ex}}(t_\star,\cdot)-f_{\mathrm{num}}(t_\star,\cdot)||_{L^2_v}$ indicated simply by ${\EuFrak{E}_f}$, as well as the computation times. One observes that in long-time regimes, the $RC$- and the $GS$-spectral methods clearly take the overhand over standard {\it FD}-methods. At shorter times however, the most advantageous method for both $\kappa$-values remains the ${\it FD}$-scheme, due to the fact that, in the initial layer the spectral schemes require (in order to be precise) many expansion terms, leading thus to longer simulation times.

\begin{table}[htb] 
\caption{Simulation times for the {\it FD}-scheme, {\it RC}- and {\it GS}-scheme (impl. Euler discretization), at instant $t_\star=10$ and several sets of parameters. The discretization step is fixed at $\Delta t=0.01$.}
\tiny
\begin{tabular}{c||c|c|c|c||c|c|c|c|}
\!\!\!$t_\star=10$&Param.&${\EuFrak{E}_m}$&${\EuFrak{E}_f}$& time&Param.& ${\EuFrak{E}_m}$&${\EuFrak{E}_f}$& time\!\!\!\\
\hline
\!\!\!$FD$($\kappa=3$)\!\!\!&$N_v=100$&$1.78\cdot10^{-4}$&$6.98\cdot10^{-3}$&49&$N_v=250$&$1.96\cdot10^{-4}$&$1.09\cdot 10^{-3}$&116\!\!\!\\
\hline
\!\!\!$RC$ ($\kappa=3$)\!\!\!&$N_{RC}=8$&$3.29\cdot 10^{-6}$&$1.15\cdot10^{-3}$&0.3&$N_{RC}=10$&$3.29\cdot 10^{-6}$&$1.57\cdot10^{-3}$&0,7\!\!\!\\
\hline
\hline
 \!\!\!$FD$($\kappa=31$)\!\!\!&$N_v=100$&$4.3\cdot10^{-14}$&$9.15\cdot10^{-3}$&26&$N_v=250$&$1.11\cdot10^{-14}$&$1.43\cdot10^{-3}$&341\!\!\!\\
\hline
\!\!\!$GS$ ($\kappa=31$)\!\!\!&$N_{GS}=6$&$1.11\cdot10^{-15}$&$2.38\cdot10^{-4}$&0.17&$N_{GS}=8$&$1.11\cdot10^{-15}$&$2.38\cdot10^{-4}$&0.21\!\!\!\\
\hline
\end{tabular}
\normalsize
 \label{tab_t10}
\end{table}

\begin{table}[htb] 
\caption{Simulation times for the {\it FD}-scheme, {\it RC}- and {\it GS}-scheme (impl. Euler discretization), at instant $t_\star=10$ and several sets of parameters. The discretization step is fixed at $\Delta t=0.01$.}
\tiny
\begin{tabular}{c||c|c|c|c||c|c|c|c|}
\!\!\!$t_\star=10$&Param.&${\EuFrak{E}_m}$&${\EuFrak{E}_f}$& time&Param.& ${\EuFrak{E}_m}$&${\EuFrak{E}_f}$& time\!\!\!\\
\hline
\!\!\!$FD$($\kappa=3$)\!\!\!
&$N_v=100$&$1.12\cdot10^{-4}$&$6.96\cdot10^{-3}$&23.7
&$N_v=250$&$1.26\cdot10^{-4}$&$1.07\cdot10^{-3}$&104.9\!\!\!\\
\hline
\!\!\!$RC$ ($\kappa=3$)\!\!\!
&$N_{RC}=8$&$9.35\cdot10^{-11}$&$4.65\cdot10^{-7}$&0.094
&$N_{RC}=10$&$9.35\cdot10^{-11}$&$5.14\cdot10^{-7}$&0.109\!\!\!\\
\hline
\hline
\!\!\!$FD$($\kappa=31$)\!\!\!
&$N_v=100$&$4.16\cdot10^{-14}$&$9.17\cdot10^{-3}$&19.7
&$N_v=250$&$6.66\cdot10^{-15}$&$1.43\cdot10^{-3}$&87.0\!\!\!\\
\hline
\!\!\!$GS$ ($\kappa=31$)\!\!\!
&$N_{GS}=6$&$4.44\cdot10^{-16}$&$4.32\cdot10^{-10}$&0.203
&$N_{GS}=8$&$3.33\cdot10^{-16}$&$4.87\cdot10^{-10}$&0.313\!\!\!\\
\hline
\end{tabular}
\normalsize
\label{tab_t10}
\end{table}

%%%%%%%%%%%%%%%
Let us finish our comparison studies, by investigating some physical and mathematical properties of our spectral solutions. Firstly, let us remark that the mass computation  is very precise (machine accuracy) with the {\it FD}-scheme for large $\kappa$-values, which comes from the fact that we have a more concentrated distribution function (close to a Maxwellian), such that no mass is lost when truncating the domain. Compare this with the conservation of mass property of the {\it FD}-scheme for small $\kappa=3$ (see Fig. \ref{INIT}). The spectral schemes preserve the mass equally accurately for all sets of parameters, a property that comes from our mass-computation procedure, synthesized here. The mass is computed by integrating the spectral expansions over the velocity space, {\it i.e.} ${\mathfrak n}=\langle f \rangle$, with
$$
\langle f_{RC} \rangle =\sum_{k=0}^{{\tilde N_{RC}}} {\mathfrak a}_{2k}(t)\, \langle \Theta_{2k} \rangle \,, \qquad 
\langle f_{GS}\rangle =\sum_{k=0}^{N_{GS}} \alpha_k(t)\, \langle p_k\, f_{\kappa,a}\rangle\,,
$$
where we recall firstly that 
$$
\langle p_k\, f_{\kappa,a}\rangle= (p_k,p_0)_{f_{\kappa,a}}= \gamma_0\, \delta_{k,0} \quad \Rightarrow \quad \langle f_{GS}\rangle = \alpha_0 (t)\,  \gamma_0\,\,.
$$
and secondly
$$
\Theta_{2k} = {\mathcal C}_{2k}\, \Upsilon_\kappa = {\mathcal C}_{2k}\, { (c_\kappa/\sqrt{2\, \kappa})^{1/2}\over (1+{v^2 \over 2\, \kappa})^{\kappa+1 \over 2}}=(c_\kappa\, \sqrt{2\, \kappa})^{1/2}\, {\mathcal C}_{2k}\, {1\over (1+{v^2 \over 2\, \kappa})^{\kappa-1 \over 2}}\, \sigma_\kappa\,,
$$
which, for $\kappa=3$ and with \eqref{expi}, yields 
$$
%\begin{array}{lll}
\ds \langle f_{RC}  \rangle =\ds  \sqrt{c_\kappa\, \sqrt{2\, \kappa}}\,  \sum_{k=0}^{{\tilde N_{RC}}} { {\mathfrak a}_{2k}(t) \over 2} \, \left[ ({\mathcal C}_{2k},{\mathcal C}_{0})_{\sigma_\kappa} - ({\mathcal C}_{2k},{\mathcal C}_{2})_{\sigma_\kappa} \right] = \sqrt{c_\kappa\, \sqrt{2\, \kappa}}\, {\pi \over 2} \left( {\mathfrak a}_{0}- {{\mathfrak a}_{2} \over 2}\right)\,.
%\end{array}
$$

\begin{figure}[ht]
\begin{center}
  \includegraphics[height=4cm,width=0.48\textwidth]{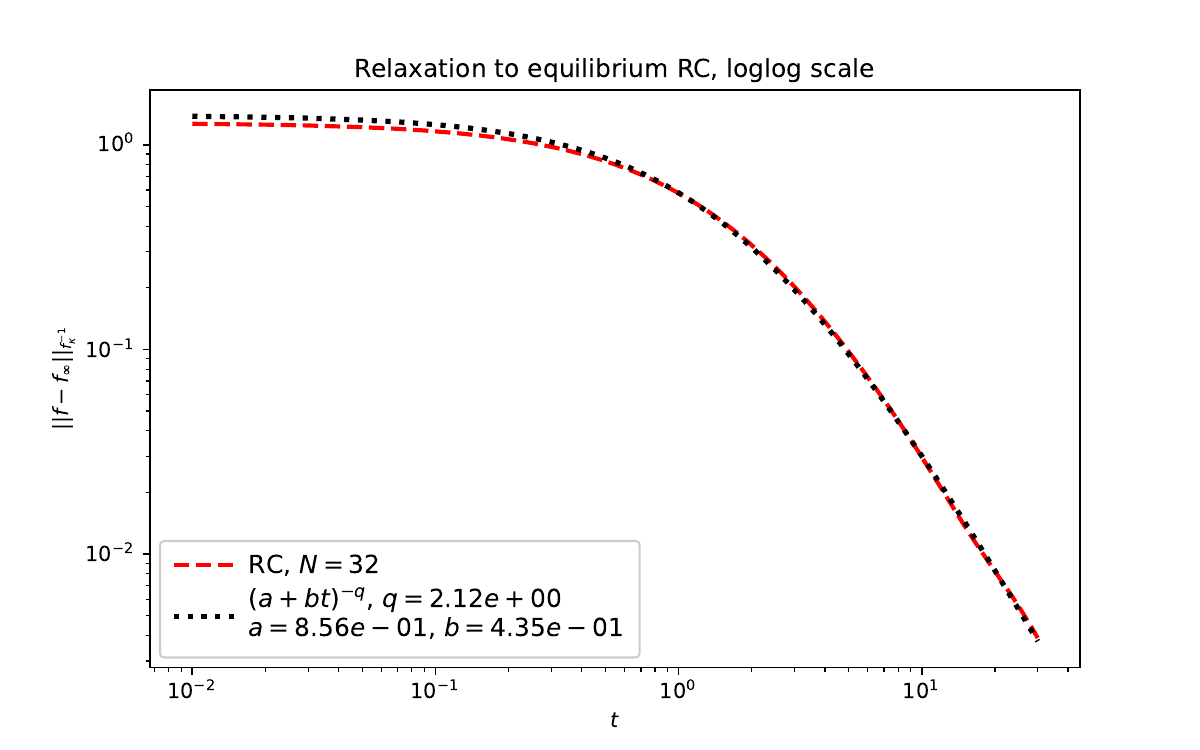}\hfill
  \includegraphics[height=4cm,width=0.48\textwidth]{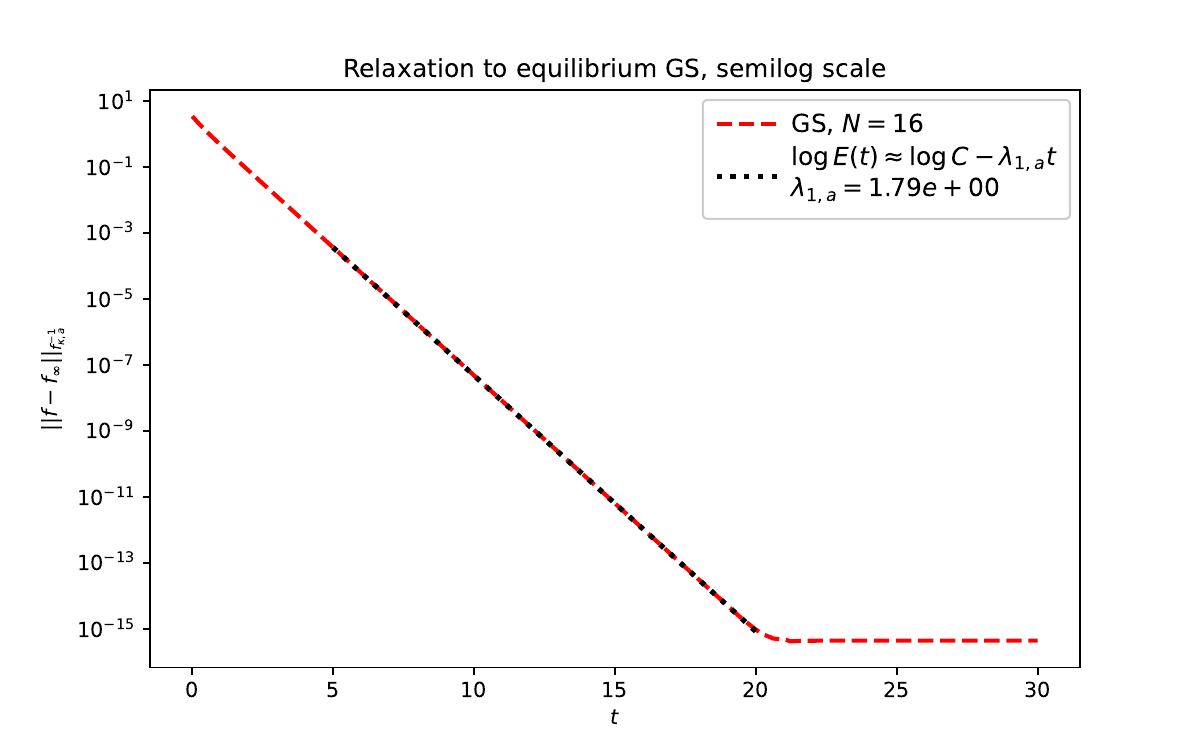}
\end{center}
\caption{Relaxation of $f(t,v)$ towards the equilibrium $f_{\infty}(v)$ as $t \rightarrow \infty$. Error plots $||f(t)-f_{\infty}||_{L^2_{f^{-1}_\kappa}}$ for the $RC$-scheme ($\kappa=3$) on the left, and the $GS$-scheme ($\kappa=31$) on the right. Superposed are the fitting curves $g(t)=(a+b\, t)^{-q}$ (left) and $E(t)=c-\lambda\, t$ (right).}\label{asymp}
\end{figure}

Secondly, the convergence of the solution towards the equilibrium is illustrated in Fig. \ref{asymp} for both schemes. The $log-log$ scale curve corresponding to the error in the $RC$-scheme is matching again very well with a function of the type $g(t)=(a+b\, t)^{-q}$, whereas the error curve corresponding to the $GS$-scheme is matching very well with a function of the type $E(t)=c-\lambda\, t$, permitting to underline the algebraic ($RC$) or exponential ($GS$) time-decay, however the decay rates are better than expected from the  mathematical Theorems \ref{Theorem_Dolbeault} and \ref{Thm_kappa}, as explained earlier in Section \ref{SEC51}.

\begin{figure}[ht]
\begin{center}
  \includegraphics[scale=0.4]{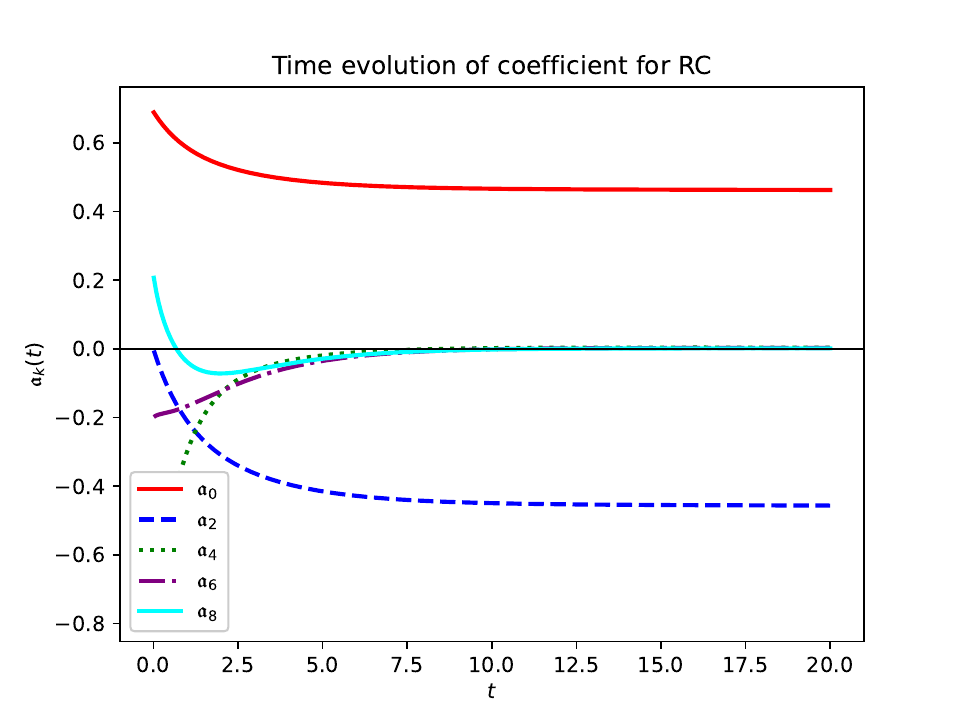}\hfill
  \includegraphics[scale=0.4]{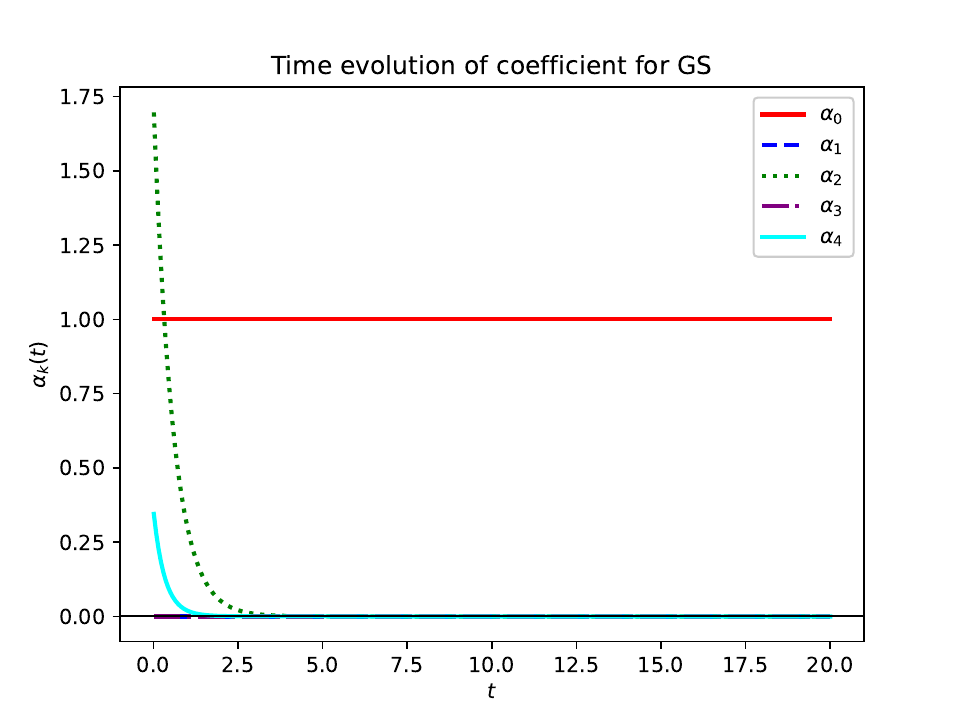}
\end{center}
\caption{Time-evolution of the first spectral coefficients ${\mathfrak a}_{k}(t)$ as well as} $\alpha_k (t)$. Left: RC-scheme, $\kappa=3$. Right: GS-scheme, $\kappa=31$. Furthermore, we took $N_{RC}=N_{GS}=10$.\label{coef_time}
\end{figure}

\begin{figure}[ht]
\begin{center}
  \includegraphics[height=4cm,width=0.48\textwidth]{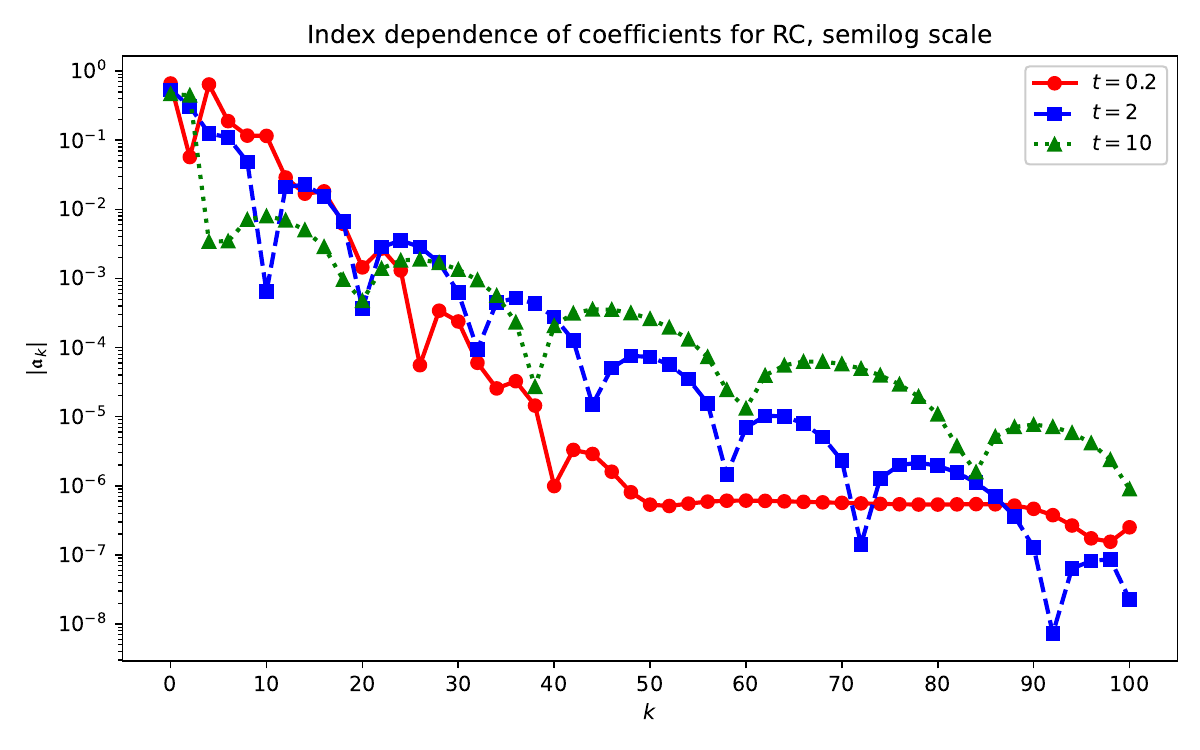}\hfill
  \includegraphics[height=4cm,width=0.48\textwidth]{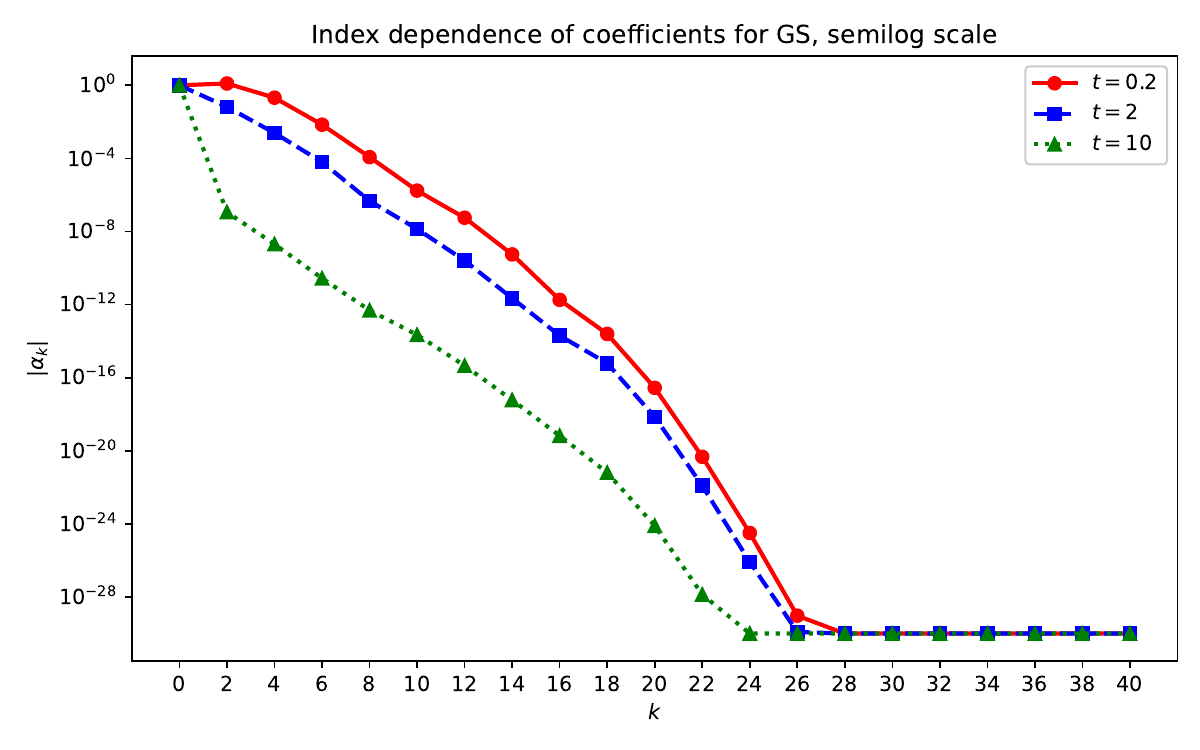}
\end{center}
\caption{Decay of the spectral coefficients ${\mathfrak a}_{k}(t_\star)$ and $\alpha_k (t_\star)$ wrt. the $k$-index, at instants $t_\star=0.2;\,\, 2;\,\, 10$. Left: RC-scheme, $\kappa=3$ and $N_{RC}=20$. Right: GS-scheme, $\kappa=31$ and $N_{GS}=10$.}\label{coef_index}
\end{figure}

Interesting is also to investigate how the spectral coefficients, namely ${\mathfrak a}_{2k}(t)$ for the {\it RC}-method and $\alpha_k(t)$ for the {\it GS}-scheme, behave in time as well as with respect to the index $k$. The time-dependence of the first coefficients is plotted on Fig. \ref{coef_time}. One observes that for the $RC$-scheme only two coefficients survive in the asymptotic limit $t \rightarrow \infty$, namely  ${\mathfrak a}_{0}(t) \rightarrow_{t \rightarrow \infty} 0.4626$ and ${\mathfrak a}_{2}(t)\rightarrow_{t \rightarrow \infty} -0.4567$, which is consistent with the fact, that we took $\kappa=3$. The $GS$-scheme requires in the asymptotic limit $t \rightarrow \infty$ only the first coefficient $\alpha_0(t)\rightarrow_{t \rightarrow \infty}1.000525$. To reduce complexity, it could be a good idea to adapt the number of used coefficients as time goes on, meaning use a large number of expansion terms in the initial time layer, whereas use fewer and fewer terms as the solution approaches the equilibrium. Such an adaptive strategy has been used in \cite{FC} by one of the authors. The magnitude of the different coefficients at a fixed instant $t_\star=0.2; \,\, 2;\,\, 10$ is illustrated in Fig. \ref{coef_index} in semi-log scale. These plots permit us to see how  many coefficients are needed at different time instants in the simulation, as well as to evaluate the decay of the coefficients with respect to the index $k$, the rate seeming to be exponential, with the particularity that the $RC$-coefficients oscillate additionally.\\

To conclude, the {\it FD}-scheme is relatively easy to implement, but inferior in accuracy and simulation time in the long-time regime, as compared to our spectral methods.  The $RC$-scheme is very efficient (accurate + rapid) for small $\kappa$-values and long-time simulations, while the $GS$-scheme is very efficient for large $\kappa$-values and long-time simulations. In the limit $\kappa \rightarrow \infty$ the {\it GS}-scheme converts into the well-known Hermite spectral scheme, used for standard Fokker-Planck equations with Maxwellian equilibria.

%%%%%%%%%%%%%%%%%%%%%%%%%%%%%%%%%%%%%%%%%%%%%%%%%%%%%%%%%%%%%%%%%%%
\section{Concluding remarks and perspectives} \label{SECC}
%%%%%%%%%%%%%%%%%%%%%%%%%%%%%%%%%%%%%%%%%%%%%%%
Let us conclude this paper by summarizing what was achieved in this work and what remains still to be done in future works. The main part of this work was concerned with the design of two efficient spectral schemes for the resolution  of a specific Fokker-Planck equation, whose stationary states are given by $\kappa$-distributions, which are (thermal) non-equilibrium distributions and describe the energetic particle population in a fusion plasma gas or in space plasmas. At the heart of spectral methods is the fact that any smooth function $f(v)$ can be expanded in the form
$$
f(v)=\sum_{k=0}^\infty \alpha_k\, \varphi_k(v)\,, \quad \forall v \in \RR\,,
$$
where $\{ \varphi_k \}_{k \in \NN}$ are global basis functions. An efficient scheme (accurate, stable and rapid) is achieved when the basis functions are adequately chosen.  We presented in this paper two different polynomial basis sets for the resolution of the energetic Fokker-Planck equation, and compare these schemes with a standard finite-difference scheme. Our conclusion is that spectral schemes (if well-designed) are superior in performance than local methods (FD), when long-time simulations are performed. For short-time simulations, standard {\it FD}-schemes are unbeatable.
However, exponential accuracy is achieved for spectral schemes only when the basis functions are adapted to the specific problem, lacking such imperative can have a drastic impact on the convergence rate. In a future work we shall come closer to the physical reality, by including the diffusion coefficient $D(v)$ in the Fokker-Planck operator, as well as introducing a transport term (coupling with Poisson equation). Addressing the transport part is not so difficult, splitting techniques permit to separate the transport part from the collisional part treated in the present paper. However, considering a velocity-dependent diffusion coefficient $D(v)$ is more intricate and the here presented spectral schemes have to be adapted to the new, more physical problem.

\bigskip

%\newpage

	%%%%%%%%%%%%%%%%%%%%%%%%%%%%%

\end{document}

%% file: ex_shared.tex
\usepackage{lipsum}
\usepackage{amsfonts}
\usepackage{graphicx}
\usepackage{epstopdf}
\usepackage{algorithmic}
\ifpdf
  \DeclareGraphicsExtensions{.eps,.pdf,.png,.jpg}
\else
  \DeclareGraphicsExtensions{.eps}
\fi

\newsiamremark{remark}{Remark}
\newsiamremark{hypothesis}{Hypothesis}
\crefname{hypothesis}{Hypothesis}{Hypotheses}
\newsiamthm{claim}{Claim}
\newsiamremark{fact}{Fact}
\crefname{fact}{Fact}{Facts}

\headers{Spectral scheme for an energetic Fokker-Planck equation}{C. Negulescu, H. Parada}

\title{Spectral scheme for an energetic Fokker-Planck equation with  $\kappa$-distribution steady states
  \thanks{Submitted to the editors DATE.
\funding{This work has been carried out within the framework of the EUROfusion Consortium, funded by the European Union via the Euratom Research and Training Programme (Grant Agreement No 101052200 — EUROfusion). Views and opinions expressed are however those of the author(s) only and do not necessarily reflect those of the European Union or the European Commission. Neither the European Union nor the European Commission can be held responsible for them. H. Parada was supported by the Agence Nationale de la Recherche, Labex CIMI under grant agreement ANR-11-LABX-0040 and is currently funded by QuBiCCS project ANR-24-CE40-3008.}}}

\author{Claudia Negulescu\thanks{Universit\'e de Toulouse \& CNRS, UPS, Institut de Math\'ematiques de Toulouse UMR 5219, F-31062 Toulouse, France
  (\email{claudia.negulescu@math.univ-toulouse.fr})} \and Hugo Parada\thanks{Universit\'e de Toulouse \& CNRS, UPS, Institut de Math\'ematiques de Toulouse UMR 5219, F-31062 Toulouse, France
  (\email{hugo.parada@math.univ-toulouse.fr})}}

\usepackage{amsopn}